\documentclass[review,hidelinks,onefignum,onetabnum]{siamart220329}

\usepackage{amsmath,amsfonts}
\usepackage{amssymb}
\usepackage{amsmath}
\usepackage{algorithmic}
\usepackage{algorithm}
\usepackage{array}
\allowdisplaybreaks[4]
\usepackage{textcomp}
\usepackage{subfigure}
\usepackage{stfloats}
\usepackage[justification=centering]{caption}
\usepackage{mathrsfs}
\usepackage{url}
\usepackage{verbatim}
\usepackage{graphicx}
\usepackage{epsfig}
\usepackage{cite}
\usepackage{hyperref}
\usepackage{overpic}
\usepackage{lipsum}
\usepackage{epstopdf}

\numberwithin{equation}{section}
\newtheorem{thm}{Theorem}
\renewcommand{\theequation}{\arabic{section}.\arabic{equation}}

\newtheorem{asum}{Assumption}
\newtheorem{lem}{Lemma}
\newtheorem{cor}{Corollary}
\newtheorem{rem}{Remark}
\newtheorem{exmp}{Example}
\newtheorem{defn}{Definition}
\allowdisplaybreaks[3]

\nolinenumbers

\headers{Distributed  Non-convex Optimization}{Qinlong Lin, Yang Liu, Jianquan Lu, and Weihua Gui}

\title{Predefined-Time Distributed  Non-convex Optimization via a Time-Base Generator\thanks{Submitted to the editors DATE.
		\funding{This work was partially supported by the National Natural Science Foundation of China under grant 62173308, and the Jinhua Science and Technology Project under grant 2022-1-042.}}}

\author{Qinlong Lin\thanks{ School of Mathematical Sciences, Zhejiang Normal University, Jinhua, 321004, China 
		(\email{linqinlong@zjnu.edu.cn}).}
	\and Yang Liu\thanks{Corresponding author. Key Laboratory of Intelligent Education Technology and Application of Zhejiang Province, the School of Mathematical Sciences, Zhejiang Normal University, Jinhua, 321004, China 
		(\email{liuyang@zjnu.edu.cn}).}
	\and  Jianquan Lu\thanks{School of Mathematics, Southeast University, Nanjing, 210096, China
		(\email{jqluma@seu.edu.cn}).}
	\and  Weihua Gui\thanks{School of Automation, Central South University, Changsha, 410083, China
		(\email{gwh@csu.edu.cn}).}
}

\usepackage{amsopn}

\begin{document}
	
	\maketitle
	
	% REQUIRED
	\begin{abstract}
		In this paper, we propose two novel multi-agent systems for the resource allocation problems (RAPs) and consensus-based distributed optimization problems.
		%%%%%%%%%%%%%%%%%%%5
		Different from existing distributed optimal approaches, we propose the new time-base generators (TBGs) for predefined-time non-convex optimization.
		%%%%%%%%%%5 
		Leveraging the proposed time-base generator, we study the roughness and boundedness of Lyapunov function based on TBGs. 
		%%%%%%%%%%%%%%%%%%%%%%55
		We prove that our approach achieves predefined-time approximate convergence to the optimal solution if the cost functions exhibit non-strongly convex or even non-convex characteristics.
		%%%%%%%%%%%%%%%55
		Furthermore, we prove that our approaches converge to the optimal solution if cost functions are generalized smoothness, and exhibit faster convergence rate and CPU speed.
		%%%%%%%%%%%%%%%%%%%%%%%%%%%%%%
		Finally, we present numerous numerical simulation examples to confirm the effectiveness of our approaches.
	\end{abstract}

	
	% REQUIRED
	
	\section{Introduction}
	%第一段分布式优化
	In recent years, the rapid advancements in large-scale systems and networks have sparked significant interest in distributed optimization problems.
	The distributed optimization method, as an effective optimization technique, exhibits several desirable features, such as, reduce communication consumption, strong fault-tolerance, safer operation, and greater privacy.
	Consequently, distributed optimization finds applications in various practical problems, such as, distributed resource allocation problems (RAPs) \cite{ZD-2018,SL-2018}, multiobjective distributed optimization problems \cite{YL-2023,KZ-2023}, and economic dispatch problems \cite{CL-2018}. Additional applications can be found in the works of \cite{XG-2020,AN-2009,TY-2018}.
	
	%	%第二段资源分配问题
	The RAPs and the consensus-based distributed optimization are two common problems in the field of optimization.
	The RAPs stand as a pivotal optimization problem across various applications, such as economic dispatch \cite{SS-2023}, power systems \cite{PY-2015} and communication networks \cite{DM-2016}.
	To address the RAPs, a large number of methods have been introduced, such as genetic algorithm \cite{NH-2009}, modified particle swarm algorithm \cite{HY-2017}, Pareto front method \cite{JQ-2019}. 
	As the scale of networks expands, fully distributed allocation algorithms have emerged as a highly favorable alternative, offering a resilient, scalable, and privacy-preserving approach to address these pressing issues.
	For instance, distributed Newton algorithms and accelerated dual descent algorithms, both exhibiting superlinear convergence rates, have been developed by researchers \cite{EA-2018,MA-2018}.
	Additionally, an initialization-free algorithm has been proposed for addressing constrained RAPs \cite{PY-2018}, while a game-based approach has been explored for resource allocation with elastic supply \cite{RS-2018}. 
	Consensus-based distributed optimization problem has gained significant attention in recent years.  
	Historically, consensus-based distributed optimization has been applied in various contexts, such as networked vehicle/UAV coordination/control, distributed estimation and learning, as well as information processing in sensor networks.
	Some examples include the multiobjective distributed optimization \cite{YL-2023}, distributed averaging \cite{KH-2014}, power system control \cite{LU-2013}.  
	In particular, in order to get a set of weak Pareto optimal solutions, multiobjective optimization problems can be converted into their equivalent distributed optimization problems \cite{YL-2023}. 
	
	%3多智能体
	
	%
	Multi-agent systems (MASs) have garnered significant attention from scholars, leading to the emergence of numerous distributed optimization methods based on MASs in the last two decades. Some methods are diverse, with some addressing constrained distributed optimization problems, while others focusing on the types of communication graphs, as detailed in the following references \cite{QJ-2015,XP-2017,SJQ-2019,SXY-2019,HXY-2019,IA-2019,BJ-2019,DSJ-2019}. 
	In a wide range of practical optimal control problems, agents are tasked with swiftly reacting to optimal solutions and achieving the desired outcomes of the corresponding tasks at a prescribed time limit.
	For instance, given the intermittency and uncertainty inherent in renewable resources within power systems, it becomes imperative to devise a swift economic dispatch algorithm that guarantees efficient and dependable operation.
	Therefore, the convergence speed of distributed optimization methods deserves further investigation.
	Numerous scholars have conducted research on finite/fixed-time stability and extended it to various systems, such as time-delay systems, discontinuous homogeneous systems, uncertain systems, and so forth, one can consult \cite{PED-2019,CURH-2020,Liu-2014,YW-2012,zhudu-2024}.
	In recent years, a large number of approaches have been introduced for finite/fixed-time distributed optimization \cite{HL-2023,PW-2017,GZ-2018,HWW-2022,PWJ-2017,MA-2022,BQZA-2019}.
	Existing works on finite/fixed-time optimization methods have several limitations: (1) The convergence time is dependent on the initial conditions; (2) The actual required time interval is shorter than the estimated one; (3) The relationship between the protocol parameters and the convergence time is not straightforward, making it difficult to directly calculate a predefined time. These limitations are also discussed in \cite{YL-2023}.
	Specifically, the settling time in \cite{YL-2023,KZ-2023,PSC-2024,LIM-2023,CYURBH-2020,GC-2022,GZ-2020,LS-2022}, which relies on the predefined/prescribed-time convergence approaches, can be set with arbitrary precision.
	From these papers, we can observe that these algorithms indeed demonstrate excellent convergence effects.
	
	%%%%%%%%%动机
	However, upon reviewing articles that investigate predefined-time algorithms, it becomes evident that strong convexity is an indispensable property within the stipulated conditions \cite{YL-2023,KZ-2023,GZ-2020,LS-2022,GC-2022}. 
	If cost functions are not strongly convex, the MASs proposed in previous works are unlikely to converge to the optimal solution, and they cannot converge to a local optimal solution under the case that the optimal solution is not unique. 
	Furthermore, there is no MASs proposed yet that can prove convergence to the optimal solution within a predefined time if cost function are not strongly convex. Therefore, it is worth exploring how to relax this condition.
	According to \cite{SL-2018,GZ-2020}, we recognize that a standard MAS designed to solve predefined optimization problems can be categorized into two components: a system of differential equations aimed at identifying equilibrium points, and a TBG utilized to expedite the convergence rate. 
	However, although these two parts can be combined to achieve predefined-Time Distributed Optimization, they are still independent of each other in some ways. In other words, the role of TBG is merely to improve efficiency, and it has little impact on the convergence properties of the differential equation dynamic system.
	
	In this paper, we propose the novel system of differential equations and the new version of TBGs. Combing the above two parts, two new MASs are introduce to solve RAPs and consensus-based distributed optimization problem, respectively.
	Drawing upon the preceding discussion, the key contributions of this paper can be concluded as follows:
	
	1) A new type of stability called $\mu$ contraction is proposed. After that, we introduce the new TBGs, which are different from \cite{KZ-2023,GZ-2020,LS-2022,BQZA-2019,GC-2022}.
	With the new TBG, we propose a new type of Lyapunov theory different from \cite{YL-2023,KZ-2023,GZ-2020,LS-2022,GC-2022}. Advantageously over the TBG mentioned in the previous papers, the robustness of Lyapunov stablility is derived under the new TBG. Furthermore, we also demonstrate that the Lyapunov function, subject to bounded nonlinear disturbances, can still converge to a controllable and bounded neighborhood.
	
	2) Leveraging the roughness of the Lyapunov function, we design new predefined-time MASs for RAPs and consensus-based distributed optimization problems. We prove that the proposed MASs converge to an optimal solution at the predefined time. However, different from \cite{YL-2023,KZ-2023,GZ-2020,LS-2022,GC-2022}, the condition that local cost functions are strong convexity is no longer a necessity, as we have generalized our conclusion to the case of strictly convex, convex or even non-convex.
	
	3) Under the boundedness of Lyapunov function, the proposed approach is proven to be convergent to an optimal solution at the predefined time if the cost functions are generalized smoothness, and exhibit faster convergence rate and faster CPU speed.
	
	4) Some examples are presented for a comprehensive comparison with existing works. In contrast to centralized optimization methods, the proposed MASs maintain the benefits of distributed optimization algorithms. Additionally, when compared to \cite{GZ-2020,YL-2023,LS-2022,MA-2022}, our approach offers a simpler mechanism for setting the predefined time. Furthermore, by integrating it with standard predefined-time MASs \cite{YL-2023,KZ-2023,GZ-2020,LS-2022,GC-2022}, we effectively alleviate certain necessary conditions, enhancing its practical applicability and flexibility.

	In Section \ref{s2}, we introduce some notations and preliminaries.
	In Section \ref{s3}, we introduce two types of optimization problems and define the predefined-time optimization. 
	In Section \ref{s4}, We propose several theorems related to TBG.  Furthermore, we introduce two predefined-time MASs tailored for distributed optimization problems and demonstrate their stability.
	In Section \ref{s5}, we provide numerical examples to illustrate the effectiveness of our proposed approaches.
	In Section \ref{s6}, we conclude this paper by summarizing our main findings and contributions.
	
	\section{Preliminaries}\label{s2}
	
	{\it Notations}: The sets $\mathbb{R}$, $\mathbb{R}^{n}$ and $\mathbb{R}^{+}$ represent the collection of real numbers, the set of $n$-dimensional real vectors and the
	set of positive numbers, respectively. 
	$\parallel \cdot\parallel$ denotes the Euclidean norm. 
	${\rm col}[z_{1}, \cdots, z_{n}]:= [z_{1}^{T} ,\cdots, z_{n}^{T} ]^{T}$.
	$\otimes$ and $\times$ denotes the Kronecker product operator and the Cartesian product operator, respectively. 
	{\it Graph theory fundamentals}: Let $\mathcal{G} = (\mathcal{V},\mathcal{E},B)$ denote a graph. $\mathcal{V} =P_{N}:= \{1, \cdots, N\}$ is the set of nodes in graph $\mathcal{G}$. $\mathcal{E}\in \mathcal{V}\times \mathcal{V}$ is the edges set in graph $\mathcal{G}$. $B =[b_{ij}]\in \mathbb{R}^{N\times N} $ is a weighted adjacency matrix. $b_{ij}>0$ if $(i, j) \in \mathcal{E}$, and $b_{ij}=0$ otherwise. Let $N_{i} = \{j \in V: b_{ij} \ne 0\}$ denote the
	set of the neighbors of $i$. $L$ denotes the Laplacian matrix of the graph $\mathcal{G}$ with $L_{ii} = \sum_{j=1}^{N} b_{ij}$ ($\forall i\in P_{N}$). $L_{ij}=-b_{ij}$ for $i \ne j$. 
	
	For a continuously differentiable function $f$: $\mathbb{R}^{n} \to \mathbb{R}$, we introduce several definitions as follows.
	
	\begin{defn}(see \cite{HXY-2019})
		For $\forall z_{1}$, $z_{2}$, $f$ is $M$-smooth
		if  
		$
		\parallel \nabla f(z_{1})-\nabla f(z_{1})\parallel \le M\parallel z_{1}-z_{2}\parallel,
		$
		where $M > 0$.
	\end{defn}
	\begin{defn}(Generalized smoothness)
		For $\forall z_{1}$, $z_{2}$, $f$ is generalized smooth if 
		$
		\parallel \nabla f(z_{1})- \nabla f(z_{2})\parallel \le M\parallel z_{1}-z_{2}\parallel+\widetilde{M},
		$
		where $M$, $\widetilde{M} \ge 0$ and $M+ \widetilde{M} \ne 0$.
	\end{defn}
	
	\begin{defn}(see \cite{RAP-2006})
		For $\forall z_{1}$, $z_{2}$, $f$ is convex if  $f(\theta z_{1}+(1-\theta)z_{2})\le \theta f(z_{1})+(1-\theta)f(z_{2})$ for $\forall\theta\in (0,1)$.
	\end{defn}
	
	\begin{defn}(see \cite{RAP-2006})
		For $\forall z_{1}\ne z_{2}$, $f$ is strictly convex if  $f(\theta z_{1}+(1-\theta)z_{2})< \theta f( z_{1})+(1-\theta)f(z_{2})$ for $\forall \theta\in (0,1)$.
	\end{defn}
	
	\begin{defn}(see \cite{HXY-2019})
		For $\forall z_{1}$, $z_{2}$,	$f$ is said to be strongly convex if $f(z_{2})>  f(z_{1})+\nabla f(z_{1})^{T}(z_{2}-z_{1})+\frac{l}{2}\parallel z_{2}-z_{1}\parallel^{2}$, where $l >0$.
	\end{defn}	
	
	\section{Problem formulation}\label{s3}
	\subsection{Optimization Problems}
	In this paper, we study the two type of optimization problems: RAPs and consensus-based distributed optimization problems, i.e.
	\begin{equation}\label{xd}
		\begin{split}
			\min f(x)=\sum_{i=1}^{N}f_{i}(x_{i}),\;\;\;
			\textbf{s.t.} \;\sum_{i=1}^{N}x_{i}=\sum_{i=1}^{N}q_{i}=q_{0}.
		\end{split}
	\end{equation}
	and 
	\begin{equation}\label{xx}
		\begin{split}
			\min f(x)=\sum_{i=1}^{N}f_{i}(x_{i}),\;\;\;
			\textbf{s.t.} \;x_{i}=x_{j},
		\end{split}
	\end{equation}
	where $x_{i} \in \mathbb{R}^{n}$, $f_{i}(x_{i})$: $\mathbb{R}^{n} \to \mathbb{R}$; $f_{i}$ has a minimizer; $x= {\rm col}[x_{1}, x_{2},..., x_{N}]$; $q_{i}$, $q_{0} \in \mathbb{R}$.
	
	No we give the first assumption:
	\begin{asum}\label{Ass1}
		Assume that the graph between agents is undirected and connected.
	\end{asum}
	Then we have the following lemma.
	\begin{lem}\label{ll}(see \cite{YL-2023}, \cite{GZ-2020})
		Under the Assumption \ref{Ass1}, for a Laplacian matrix $L$, we have that 
		$\lambda_{N}(L)\ge\dots\ge\lambda_{2}(L)>0$, where $\lambda_{N}(L^{2})$ is the largest eigenvalue of $L^{2}$ and $\lambda_{2}(L)$ is the second smallest eigenvalue of $L$. Furthermore $x_{i} =
		x_{j}$ if and only if $Lx = 0$.
		
	\end{lem}
	
	\subsection{Predefined-Time Optimization and TBG} 
	It may be difficult to obtain accurate distributed optimization solutions because of technical difficulties and complexity or cost, while converges to the neighborhood of an optimal solution at a predefined time is a good choice, which can provide considerable benefits by offering simple and feasible designs, and potentially even enhancing performance.
	Now we first give the definitions of predefined-time approximate convergence and predefined-time optimization.
	\begin{defn}\label{xaf}(\cite{YL-2023})
		System $\dot{x}(t) = A(t,t_{p})f(x(t))$ ($\forall x(0)$ and $\forall t\ge0$)
		achieves predefined-time approximate convergence at $t_{p}$ if: \\
		1) $\lim_{t\to t_{p+}} \parallel x(t)\parallel\le \epsilon$ and $0<\epsilon=\epsilon(x(0)) \ll 1$;\\
		2) $\parallel x(t)\parallel\le \epsilon$ for $\forall t\ge t_{p}$;\\
		3) $\lim_{t\to \infty} \parallel x(t)\parallel=0$.
	\end{defn}
	
	Based on Definition \ref{xaf}, the definition of the predefined-time optimization is similar to Definition 6 in\cite{YL-2023}.
	\begin{defn}\label{xaf2}(\cite{YL-2023})
		For  $ \forall x(0)$ and $\forall t\ge0$, the continuous-time
		optimization approach achieves convergence to the optimal solution $x^{*}$ at the predefined time $t_{p}$ if \\
		1) $\lim_{t\to t_{p+}} \parallel x(t)-x^{*}\parallel\le \epsilon$ and $0<\epsilon=\epsilon(x(0)) \ll 1$;\\
		2) $\parallel x(t)-x^{*}\parallel\le \epsilon$ for $\forall t\ge t_{p}$;\\
		3) $\lim_{t\to \infty} \parallel x(t)-x^{*}\parallel=0$.
	\end{defn}
	
	Note that $A(t,t_{p})$ in Definition \ref{xaf} is a proper TBG. TBG $A(t,t_{p})$ has also been defined in different forms. In \cite{GZ-2020}, $A(t,t_{p})=v ([\dot{\gamma}(t)]/[1-\gamma(t)+\varsigma])$, where $\gamma(t)$ satisfies
	\begin{equation*}
		\begin{cases}
			\gamma(t)=0, t=0;\\
			\gamma(t)=1,t\ge t_{p},
		\end{cases}\;\;\;
		\begin{cases}
			\dot{\gamma}(t)=0, t\ge t_{p}\;\; {\rm or}\;\; t=0;\\
			\dot{\gamma}(t)>0, t_{p}>t>0.
		\end{cases}
	\end{equation*}
	Furthermore, in \cite{SL-2018}, $A(t,t_{p})=uT(t,t_{p})=u(\zeta(t,\varsigma))^{\prime} $, where $u>0$, $0<\varsigma\ll1$ and $\zeta(t,\varsigma)$ satisfies
	$
	\lim_{\varsigma\to0+}[\zeta(t_{p+},\varsigma)-\zeta(0,\varsigma)]=+\infty,
	\zeta(t,\varsigma)-\zeta(t_{p+},\varsigma)\ge0,\;\;\forall t>t_{p},
	\lim_{t\to\infty}[\zeta(t,\varsigma)-\zeta(0,\varsigma)]=+\infty.
	$
	It is noted that $\zeta(t,\varsigma)$ can be relaxed into a weaker form $\zeta(t)$ if $\lim_{\varsigma\to0+}[\zeta(t_{p+},\varsigma)-\zeta(0,\varsigma)]=a$ is sufficiently large. Take an example, $A(t,t_{p})=\zeta^{\prime}(t)$, where
	\begin{equation*}
		\zeta^{\prime}(t)=
		\begin{cases}
			99,&0\le t\le0.05,\\
			0,&t\ge0.05.
		\end{cases}
	\end{equation*}
	However, there previous TBGs can not be used to achieve predefined-time optimization if the cost functions are not strongly convex.
	
	In this paper, we propose a new TBG. Consider the following system:
	\begin{equation}\label{Vt}
		\dot{V}(t)=-A(t,t_{p})V(t),
	\end{equation}
	where $V(t)$: $\mathbb{R}\to\mathbb{R}^{n}$, $A(t,t_{p})$ is a $n\times n$ matrix.
	We denote $\Phi(t,\tau)$ by the associated evolution operator satisfied  $V(t)=\Phi(t,\tau)V(\tau)$ for $\forall t$, $\tau\in(0,t_{p}]$ or $[t_{p},+\infty)$.
	
	Assume that positive function $\mu(t,\tau,\varsigma)$ satisfy:\\
	$\mathcal{A}1)$ $\lim_{t\to +\infty}\mu(t,\tau,\varsigma)=0$, $\mu(t,\tau,\varsigma)$ is decreasing for $t$;\\
	$\mathcal{A}2)$ $\lim_{\varsigma\to 0}\mu(t_{p}^{+},\tau,\varsigma)=0$, \\
	where $t$, $\tau$, $t_{p}$, $\varsigma\in\mathbb{R}$.
	Based on the above conditions, we introduce the predefined-time $\mu$ contraction.
	\begin{defn}\label{aptu}
		System (\ref{Vt}) is said to admit a predefined-time $\mu$ contraction
		if 
		$
		\parallel\Phi(t,\tau)\parallel\le\mu(t,\tau,\varsigma)$ for $\forall t\ge\tau\ge0$. 
		Furthermore,
		$A(t,t_{p})$ is a TBG if system (\ref{Vt}) admits predefined-time $\mu$ contraction.
	\end{defn}
	\begin{rem}
		Note that $\mathcal{A}1)$-$\mathcal{A}2)$ are important conditions for ensuring that system (\ref{Vt}) achieves predefined-time approximate convergence at $t_{p}$, we can observe similar properties from previous TBGs. Recall TBG $A(t,t_{p})= vT(t,t_{p})=v(\zeta(t,\varsigma))^{\prime}$,
		we have that $\Phi(t,\tau)=e^{-v(\zeta(t,\varsigma)-\zeta(\tau,\varsigma))}$.
		Then we set $\mu(t,\tau,\varsigma)=e^{-v(\zeta(t,\varsigma)-\zeta(\tau,\varsigma))}$, $\mathcal{A}1)$- $\mathcal{A}2)$ hold, and predefined-time $\mu$ contraction reduces to predefined-time exponential contraction. Additionally, review TBG $A(t,t_{p})=v ([\dot{\gamma}(t)]/[1-\gamma(t)+\varsigma])$, 
		we have that $\Phi(t,\tau)=\bigg(\frac{1+\varsigma-\gamma(t)}{1+\varsigma-\gamma(\tau)}\bigg)^{v}$. Then set $\mu(t,\tau,\varsigma)=e^{-v(-\ln(1+\varsigma-\gamma(t))+\ln(1+\varsigma-\gamma(\tau)))}$, we can observe that $\mathcal{A}1)$ and $\mathcal{A}2)$ hold. Thus, $\mathcal{A}1)$ and $\mathcal{A}2)$ are standard.
	\end{rem}
	
	\section{Main Results}\label{s4}
	\subsection{Some theorems of TBG}
	Now, we first propose a theorem related to TBG as follows.
	\begin{thm}\label{pac}
		System (\ref{Vt}) under TBG in Definition \ref{aptu} achieves predefined-time approximate convergence at $t_{p}$.
	\end{thm}
	{\it Proof}: According to the definition of predefined-time $\mu$ contraction, we obtain that 
	$
	\parallel V(t)\parallel
	=\parallel\Phi(t,0)V(0)\parallel
	\le \mu(t,0,\varsigma) \parallel V(0)\parallel.
	$
	Hence, by $\lim_{t\to +\infty}\mu(t,\tau,\varsigma)=0$ we have that $\lim_{t\to \infty}\parallel V(t)\parallel=0$. Moreover, we also obtain that
	$
	\parallel V(t_{p+})\parallel \le \mu(t_{p+},0,\varsigma) \parallel V(0)\parallel.
	$
	According to $0<\varsigma\ll 1$ and $\lim_{\varsigma\to 0}\mu(t_{p}^{+},\tau,\varsigma)=0$, we have that $\parallel V(t_{p+})\parallel\le\epsilon$. Finally, because $\mu(t,\tau,\varsigma)$ is decreasing for $t$, we have that $\parallel x(t)\parallel\le \epsilon$ $(\forall t\ge t_{p})$.
	
	In fact, setting up a TBG that incorporates $\varsigma$ is not a trivial task. Therefore, we consider simplifying the setup of TBG, i.e. $\mathcal{A}1)$ and  $\mathcal{A}2)$ are relaxed into a weaker form without $\varsigma$. We add two positive constants $\alpha$ and $D$ such that positive function $\mu_{\alpha,D}(t,\tau)$ satisfy:\\
	$\mathcal{B}1)$ $\lim_{t\to +\infty}\mu_{\alpha,D}(t,\tau)=0$, $\mu_{\alpha,D}(t,\tau)$ is decreasing for $t$, $\alpha$ and increasing for $D$, where $t$, $\tau\in\mathbb{R}$;\\
	$\mathcal{B}2)$ $\alpha$ is sufficiently large such that $\lim_{t\to t_{p+}}\mu_{\alpha,D}(t,0)\le\epsilon$, where $\epsilon$ is small.\\
	$\mathcal{B}3)$ $\mu_{\alpha,D}(t,\tau)=D\mu^{\alpha}(t,\tau)$, $\mu(t,\tau)\mu(\tau,s)=\mu(t,s)$, $\frac{\partial\mu(0,t)}{\partial t}\cdot\mu^{-1}(0,t)\ge1 $,
	where $\mu(t,\tau)=\mu_{1,1}(t,\tau)$.
	
	We are going to take some examples to show that $\mathcal{B}1) - \mathcal{B}3)$ are standard.
	\begin{exmp}\label{E1}
		If $A(t,t_{p})=\alpha$ ($\alpha$ is sufficiently large), $D=2$. Then $\mu_{\alpha,D}(t,\tau)=2e^{-\alpha(t-\tau)}$, $\epsilon=2e^{-\alpha t_{p}}$ is small. If $A(t,t_{p})=\alpha\vartheta^{\prime}(t)\vartheta^{-1}(t)$, where $\vartheta(t)=\frac{2}{\pi}e^{ t}(\frac{\pi}{2}+\arctan t)$, then $\Phi(t,\tau)=\vartheta^{-\alpha}(t)\vartheta^{\alpha}(\tau)$, and $\mu_{\alpha,D}(t,\tau)=\vartheta^{-\alpha}(t)\vartheta^{\alpha}(\tau)$ ($\alpha$ is sufficiently large),
		hence $\epsilon=\vartheta^{-\alpha}(t_{p})$ is small. If $\mu_{\alpha,D}(t,\tau)=De^{-\alpha(t-\tau)}$, then $De^{-\alpha(t-\tau)}=D(e^{-(t-\tau)})^{\alpha}$, $e^{-(t-\tau)}\cdot e^{-(\tau-s)}=e^{-(t-s)}$ and $\frac{\partial\mu(0,t)}{\partial t}\cdot\mu^{-1}(0,t)=1 $. If $\mu_{\alpha,D}(t,\tau)=D\vartheta^{-\alpha}(t)\vartheta^{\alpha}(\tau)$, then $D\vartheta^{-\alpha}(t)\vartheta^{\alpha}(\tau)=D(\vartheta^{-1}(t)\vartheta(\tau))^{-\alpha}$, $\vartheta^{-1}(t)\vartheta(\tau)\vartheta^{-1}(\tau)\vartheta(s)=\vartheta^{-1}(t)\vartheta(s)$ and $\frac{\partial\mu(0,t)}{\partial t}\cdot\mu^{-1}(0,t)=1+\frac{1}{(\frac{\pi}{2}+\arctan t)(1+t^{2})}\ge1$. Hence, $\mathcal{B}1)$-$\mathcal{B}3)$ are standard.
	\end{exmp}
	
	Now we are going to study the robustness
	of predefined-time $\mu$ contraction.
	
	\begin{thm}\label{thr}
		Suppose that system (\ref{Vt})
		admits a predefined-time $\mu$ contraction with assumptions $\mathcal{B}1)$-$\mathcal{B}3)$, if $\parallel \widetilde{A}(t)\parallel \le\delta$ and $\delta\le\alpha/ D$, then system 
		\begin{equation}\label{AVBV}
			\dot{V}(t)=-A(t,t_{p})V(t)+\widetilde{A}(t)V(t),\;\;\forall t\ge 0
		\end{equation} 
		also admits a predefined-time $\mu$ contraction, i.e.,
		$
		\parallel\Gamma(t,\tau)\parallel\le \mu_{\alpha-\delta D,D}(t,\tau),(t\ge\tau),
		$
		where $\Gamma(t,\tau)$ is the evolution operator associated to system (\ref{AVBV}).
	\end{thm}
	
	We first give the generalized Gronwall inequality.
	\begin{lem}\label{LV}
		Let $v$ be a non-negative continuous function satisfying 	
		\begin{equation}\label{lV}
			v(t)\le \mu_{\alpha,D}(t,\tau) + \delta \int^{t}_{\tau} \mu_{\alpha,D}(t,s)\frac{\partial\mu(0,s)}{\partial s}\cdot\mu^{-1}(0,s)v(s)\mathrm{d}s,
		\end{equation}
		where $t\ge \tau$. If $\delta< \frac{\alpha}{D}$, then
		$
		v(t)\le \mu_{\alpha-D\delta,D}(t,\tau),
		$
		where $t\ge s$.
	\end{lem}
	
	{\it Proof:} By (\ref{lV}) and $\mathcal{B}3)$, we have that
	\begin{equation*}
		\begin{split}
			v(t)\mu^{-\alpha}(t,0)
			\le
			D\mu^{\alpha}(0,\tau) + \delta D\int^{t}_{\tau} \mu^{\alpha}(0,s)\frac{\partial\mu(0,s)}{\partial s}\cdot\mu^{-1}(0,s)v(s)\mathrm{d}s.
		\end{split}
	\end{equation*}
	Set $\widetilde{v}(t)=\mu^{\alpha}(0,t)v(t)$, then we have
	$
	\widetilde{v}(t)
	\le
	D\mu^{\alpha}(0,\tau) + \delta D \int^{t}_{\tau} \widetilde{v}(s)\frac{\partial\mu(0,s)}{\partial s}\cdot\mu^{-1}(0,s)\mathrm{d}s.
	$
	By Bellman inequality, we obtain that
	\begin{equation*}
		\begin{split}
			\widetilde{v}(t)
			\le
			&D\mu^{\alpha}(0,\tau)\cdot  e^{\delta D \int^{t}_{\tau}\frac{\partial\mu(0,s)}{\partial s}\cdot\mu^{-1}(0,s)\mathrm{d}s}
			=D\mu^{\alpha}(0,\tau)\cdot\mu(t,\tau)^{-\delta D}.
		\end{split}
	\end{equation*}
	Hence, by $\mathcal{B}3)$, we have that
	$
	v(t)\le \mu_{\alpha-D\delta,D}(t,\tau).
	$
	Then we obtain the desired results.
	
	Now the proof of Theorem \ref{thr} is given.
	
	{\it Proof:}
	Consider the space $\mathcal{B}(V)$. Set $t\ge \tau\ge0$
	and 
	$
	\mathfrak{G}=\{\Gamma:[0,+\infty)\times[0,+\infty)\to \mathcal{B}(V) \},
	$
	where $\Gamma$ is continuous and $\parallel \Gamma\parallel<\infty$.
	We set 
	$
	(\mathscr{L}\Gamma)(t,\tau)=\Phi(t,\tau)+\int^{t}_{\tau} \Phi(t,s)\widetilde{A}(s)\Gamma(s,\tau)\mathrm{d}s.
	$
	Note that
	\begin{equation*}
		\begin{split}
			\parallel(\mathscr{L}\Gamma)(t,\tau)\parallel
			\le&
			\parallel\Phi(t,\tau)\parallel+\int^{t}_{\tau} \parallel\Phi(t,s)\parallel\cdot\parallel \widetilde{A}(s)\parallel\cdot\parallel \Gamma(s,\tau)\parallel\mathrm{d}s\\
			\le& \mu_{\alpha,D}(t,\tau) 
			+ \delta \parallel \Gamma\parallel\int^{t}_{\tau} \mu_{\alpha,D}(t,s) \zeta(t)\mathrm{d}s,
		\end{split}
	\end{equation*}
	where $\parallel\Gamma\parallel=\sup\parallel\Gamma(t,\tau)\parallel$.
	Since $\alpha>0$, we have that
	$
	\parallel \mathscr{L}\Gamma \parallel\le D+\delta \frac{D}{\alpha}\parallel \Gamma\parallel.
	$
	Therefore, we have a well-defined operator $\mathscr{L}$. We find that for any $U_{1}$, $U_{2}$,
	$
	\parallel \mathscr{L}\Gamma_{1}-\mathscr{L}\Gamma_{2}\parallel\le \delta\frac{D}{\alpha}\parallel \Gamma_{1}-\Gamma_{2}\parallel.
	$	
	Since $\delta< \dfrac{\alpha}{D}$, the operator $L$ is a contarction. Hence, there exists a unique $\Gamma$ such that $\mathscr{L}\Gamma=\Gamma$, which thus satisfies the identity	
	$
	\Gamma(t,\tau)=\Phi(t,\tau)+\int^{t}_{\tau} \Phi(t,s)\widetilde{A}(s)\Gamma(s,\tau)\mathrm{d}s.
	$
	By Lemma \ref{LV}, we obtain the desired results.

	According to Theorem \ref{thr}, we present the Lyapunov theory.
	\begin{thm}\label{vcx}
		Suppose that utilizing generalied TBGs with $\mathcal{B}1)$-$\mathcal{B}3)$ in system 
		\begin{equation}\label{vax}
			\dot{V}(x(t))=-A(t,t_{p})V(x(t)),
		\end{equation} 
		and 
		$\widetilde{A}(t)$ in
		$
		\dot{V}(x(t))=-A(t,t_{p})V(x(t))+\widetilde{A}(t)V(x(t))
		$
		satisfies $\parallel \widetilde{A}(t)\parallel \le\delta$ and $\delta\le\alpha/ D$. If there exists a Lyapunov function $V (x)$ satisfying that $V (x)\ge c \parallel x\parallel^{2}$ with $ c> 0$, then the origin of system
		$\dot{x}(t) = -A(t,t_{p})f(x(t))$ achieves predefined-time approximate convergence at $t_{p}$.
	\end{thm}
	
	Now we study the bounded solutions of system 	
	\begin{equation}\label{AVFV}
		\dot{V}(t)=-A(t,t_{p})V(t)+F(t,V(t)),
	\end{equation}
	where $F$: $R\times R^{n}\to R^{n}$ are continuous functions.
	\begin{thm}\label{thbound}
		Suppose that system (\ref{AVFV}) admits a predefined-time $\mu$ contraction  with assumptions $\mathcal{B}1)$-$\mathcal{B}3)$. If $\parallel F(t,V(t))\parallel \le \widehat{M}$ ($\widehat{M}>0$), then system (\ref{AVFV}) has a bounded solution. Furthermore, $\lim_{t\to\infty}\parallel V(t)\parallel\le\frac{D\widehat{M}}{\alpha}$.
	\end{thm}
	
	{\it Proof:} Any solution is written as 
	$
	V(t)=\Phi(t,\tau)V(\tau)+\int^{t}_{\tau}\Phi(t,s)F(s,V(s))\mathrm{d}s.
	$
	Then by $\mathcal{B}1)$-$\mathcal{B}3)$, we have that 
	\begin{equation*}
		\begin{split}
			\parallel V(t)\parallel
			\le& \parallel \Phi(t,\tau)V(\tau)\parallel+\widehat{M}\int^{t}_{\tau}\parallel \Phi(t,s)\parallel \mathrm{d}s\\
			\le& \lambda_{\alpha,D}(t,\tau)\parallel V(0)\parallel+\widehat{M}\int^{t}_{\tau}\frac{\partial\mu(0,s)}{\partial s}\mu^{-1}(0,s)\mu_{\alpha,D}(t,s)\mathrm{d}s\\
			\le &\lambda_{\alpha,D}(t,\tau)\parallel V(0)\parallel+\dfrac{D\widehat{M}}{\alpha}.
		\end{split}
	\end{equation*}
	This complete the proof. 
	
	We also have the following theorem.
	\begin{thm}
		Suppose that utilize TBG with $\mathcal{B}1)$-$\mathcal{B}3)$ in system (\ref{vax})
		and $F$ in
		$
		\dot{V}(x(t))=-A(t,t_{p})V(x(t))+F(t,V(x(t)))
		$ 
		satisfy $\parallel F(t,V(t))\parallel \le \widehat{M}$ ($\widehat{M}>0$). If there exists a Lyapunov function $V (x)$ satisfying that $V (x)\ge c \parallel x\parallel^{2}$ with $ c> 0$, then the origin of system
		$\dot{x}(t) = -A(t,t_{p})f(x(t))$ achieves predefined-time approximate convergence to the bounded neighborhood of $x^{*}$ at $t_{p}$.
	\end{thm}
	
	Now we provide some examples.
	\begin{exmp}\label{ex2}
		Consider with $t_{p} = 0.05$ and $V(0)=10$. Several TBGs are designed as follows.\\
		\begin{equation*}
			A_{1}(t,t_{p})=
			\begin{cases}
				-100\vartheta^{\prime}(t)\vartheta^{-1}(t);&t\le0.05;\\
				-10; &t>0.05.
			\end{cases}\;\;\;
			A_{2}(t,t_{p})=
			\begin{cases}
				-100\vartheta^{\prime}(t)\vartheta^{-1}(t);&t\le0.05;\\
				-200; &t>0.05.
			\end{cases}
		\end{equation*}
		where $\vartheta(t)=\frac{2}{\pi}e^{ t}(\frac{\pi}{2}+\arctan t)$\\
		$1)$  Consider system (\ref{Vt}) with $A_{1}(t,t_{p})$;\\
		$2)$ Consider system (\ref{AVBV}) with $A_{1}(t,t_{p})$ and  $\widetilde{A}(t)=20$;\\
		$3)$  Consider system (\ref{AVBV}) with $A_{1}(t,t_{p})$ and $\widetilde{A}(t)=30$;\\
		$4)$  Consider system (\ref{AVFV}) with $A_{2}(t,t_{p})$  and $F(t,V(t))=40*\sin(t)$;\\
		$5)$  Consider system (\ref{AVFV}) with $A_{2}(t,t_{p})$   and $F(t,V(t))=50*\cos(t)$;\\
		$6)$  Consider system (\ref{AVFV}) with $A_{2}(t,t_{p})$   and $F(t,V(t))=60$.
		\begin{figure}[htph]
			\centering
			{\includegraphics[width=7cm,height=5cm]{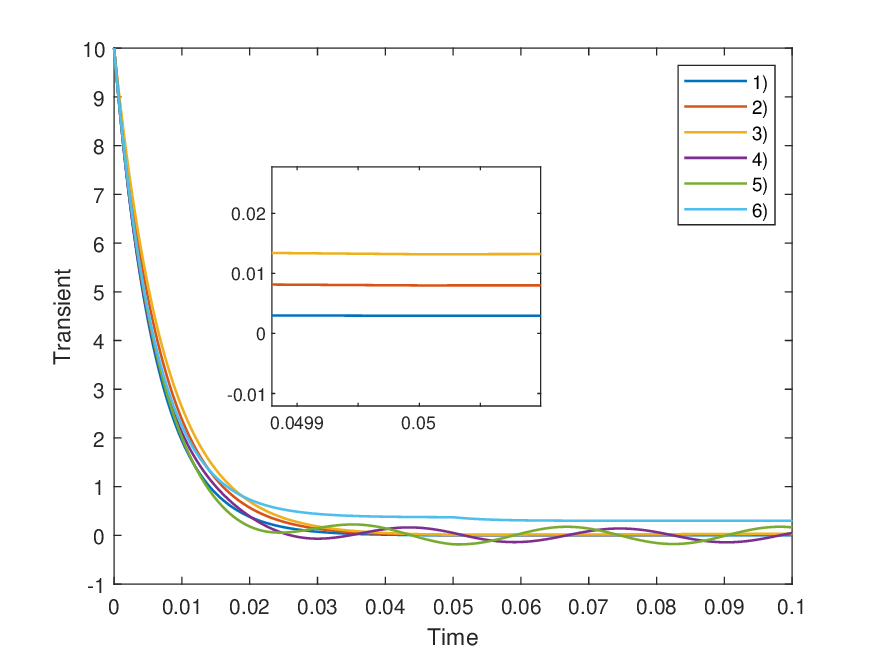}}
			\caption{Transient behaviors of $V(t)$ with different TBGs in Example \ref{ex2}.}\label{f2}
		\end{figure}
		
		Using MATLAB, Fig. \ref{f2} depicts the transient states of system (\ref{AVBV}) and system (\ref{AVFV}). According to Fig. \ref{f2}, we find that under small linear disturbances, system (\ref{AVBV}) also achieves the predefined-time approximate convergence with an exponential decline trend (see $5)$ and $6)$ in Fig. \ref{f2}). The smaller the linear disturbance received, the faster the convergence speed. Furthermore, we find that under nonlinear disturbances, system (\ref{AVFV}) achieves the predefined-time approximate convergence to a bounded neighborhood (see $1)$, $2)$ and $3)$ in Fig. \ref{f2}). Additionally, the smaller the nonlinear disturbance received, the faster the convergence speed with a narrower boundary of neighborhood.
	\end{exmp}	
	
	\subsection{Predefined-Time MAS}
	In this section, the novel predefined-time MASs under the TBG in Definition \ref{aptu} are proposed and proven. For optimization problem (\ref{xd}) and (\ref{xx}), we proposed the TBG-based predefined-Time MASs:
	\begin{equation}\label{A1}
		\begin{cases}
			\dot{x}_{i}=A(t,t_{p})\bigg(-\frac{\varrho}{A(t,t_{p})}\nabla f_{i}(x_{i})+y_{i}-x_{i}\bigg) \\
			\dot{y}_{i}=A(t,t_{p})\bigg(-\sum_{j\in N_{i}}a_{ij}(y_{i}-y_{j})+q_{i}-x_{i}
			+\frac{\varrho}{A(t,t_{p})}\sum_{j\in N_{i}}a_{ij}(z_{i}-z_{j})\bigg) \\
			\dot{z}_{i}={A(t,t_{p})}\bigg(\sum_{j\in N_{i}}a_{ij}(y_{i}-y_{j})-\sum_{j\in N_{i}}a_{ij}(z_{i}-z_{j})
			-\sum_{j\in N_{i}}a_{ij}(u_{i}-u_{j})\bigg)\\
			\dot{u}_{i}=A(t,t_{p}) \bigg(\sum_{j\in N_{i}}a_{ij}(y_{i}-y_{j})
			-\sum_{j\in N_{i}}a_{ij}(x_{i}-x_{j})\bigg),
		\end{cases}
	\end{equation}
	and 
	\begin{equation}\label{A2}
		\begin{cases}
			\dot{x}_{i}=A(t,t_{p})\bigg(-\frac{\varrho}{A(t,t_{p})}\nabla f_{i}(x_{i})-\sum_{j\in N_{i}}a_{ij}(x_{i}-x_{j})
			-\sum_{j\in N_{i}}a_{ij}(w_{i}-w_{j})\bigg) \\
			\dot{w}_{i}=A(t,t_{p})\bigg(-\sum_{j\in N_{i}}a_{ij}(x_{i}-x_{j})\bigg),
		\end{cases}
	\end{equation}
	respectively. Note that $y_{i}$, $z_{i}$, $u_{i}$, $w_{i}$ are the internal auxiliary variables; $\nabla f_{i}(x_{i})$ are the gradients of $f_{i}(x_{i})$, and $\kappa A(t,t_{p})$ (here $\forall \kappa$, $A(t,t_{p})\in\mathbb{R}^{+}$) are TBG in Definition \ref{aptu}.
	
	%For optimization problem (\ref{xd}), the first equality in (\ref{A1}) is applied to guarantee the convergence to the optimal solutions. The term $\frac{\varrho}{A(t)}$ in first equality is applied to transform the gradient $\nabla f_{i}(x_{i})$ into a sufficiently small perturbation. The terms $\sum_{j\in N_{i}}a_{ij}(x_{i}-x_{j})$, $\sum_{j\in N_{i}}a_{ij}(y_{i}-y_{j})$, $\sum_{j\in N_{i}}a_{ij}(z_{i}-z_{j})$ and $\sum_{j\in N_{i}}a_{ij}(u_{i}-u_{j})$
	%are used to realize consensus, and $q_{i}-x_{i}$ is used to guarantee the global equality constraint. 
	%Similarly, for optimization problem (\ref{xx}), the first equality in (\ref{A1}) is applied to guarantee the convergence to the optimal solutions. The term $\frac{\varrho}{A(t)}$ in first equality is applied to transform the gradient $\nabla f_{i}(x_{i})$ into a sufficiently small perturbation. The terms $\sum_{j\in N_{i}}a_{ij}(x_{i}-x_{j})$ is used to guarantee the global equality	constraint.
	
	Now we give the main results.
	\begin{thm}\label{thth1}
		Suppose that Assumption \ref{Ass1} holds. Further assume that
		$f_{i}(x_{i})$ are continuously differentiable and $M$-smooth, then MAS (\ref{A1}) converges to its equilibrium points at the predefined time $t_{p}$.
	\end{thm}
	
	{\it Proof:} 
	%Let $\nabla f(x)={\rm col}[\nabla f_{1}(x_{1}),f_{2}(x_{2}),...,f_{N}(x_{N})]$, $\nabla f(x^{*})={\rm col}[\nabla f_{1}(x_{1}^{*}),f_{2}(x_{2}^{*}),...,f_{N}(x_{N}^{*})]$, $d={\rm col}[ q_{1},q_{2},...,q_{N}]$, $y={\rm col}[ y_{1},y_{2},...,y_{N}]$,
	%$y^{*}=[ y_{1}^{*T},y_{2}^{*T},...,y_{N}^{*T}]^{T}$, $z=[ z_{1}^{T},z_{2}^{T},...,z_{N}^{T}]^{T}$, $z^{*}=[ z_{1}^{*T},z_{2}^{*T},...,z_{N}^{*T}]^{T}$, $u=[ u_{1}^{T},u_{2}^{T},...,u_{N}^{T}]^{T}$, $u^{*}=[ u_{1}^{*T},u_{2}^{*T},...u_{N}^{*T}]^{T}$. 
	Note that (\ref{A1}) can
	be rewritten as follows
	\begin{equation}\label{AA1}
		\begin{cases}
			\dot{x}=A(t,t_{p})(-\frac{\varrho}{A(t,t_{p})}\nabla f(x)+y-x) \\
			\dot{y}=A(t,t_{p})(-\mathbb{L}y+\frac{\varrho}{A(t,t_{p})}\mathbb{L}z+q-x) \\
			\dot{z}={A(t,t_{p})}(\mathbb{L}y-\mathbb{L}z-\mathbb{L}u)\\
			\dot{u}=A(t,t_{p}) (\mathbb{L}y-\mathbb{L}x),
		\end{cases}
	\end{equation}
	where $\mathbb{L}=L\otimes {\rm I}_{N} $, $q=\sum_{i\in P_{N}}q_{i}$, $\nabla f(x)={\rm col}[\nabla f_{1}(x_{1}),f_{2}(x_{2}),...,f_{N}(x_{N})]$ and $y$, $z$, $u$ are similar.
	Denote the equilibrium point of (\ref{AA1}) as $(x^{*}, y^{*},z^{*},u^{*})$. Then we obtain that 
	\begin{equation}\label{p1}
		\begin{cases}
			0=(-\frac{\varrho}{A(t,t_{p})}\nabla f(x^{*})+y^{*}-x^{*}) \\
			0=(-\mathbb{L}y^{*}+\frac{\varrho}{A(t,t_{p})}\mathbb{L}z^{*}+q-x^{*}) \\
			0=(\mathbb{L}y^{*}-\mathbb{L}z^{*}-\mathbb{L}u^{*})\\
			0=(\mathbb{L}y^{*}-\mathbb{L}x^{*}).
		\end{cases}
	\end{equation}
	Let $\hat{x}=x-x^{*}$, $\hat{y}=y-y^{*}$, $\hat{z}=z-z^{*}$ and $\hat{u}=u-u^{*}$, then we have that
	\begin{equation}\label{AAb1}
		\begin{cases}
			\dot{\hat{x}}=A(t,t_{p})(-\frac{\varrho}{A(t,t_{p})}(\nabla f(x)-\nabla f(x^{*}))+\hat{y}-\hat{x}) \\
			\dot{\hat{y}}=A(t,t_{p})(-\mathbb{L}\hat{y}+\frac{\varrho}{A(t,t_{p})}\mathbb{L}\hat{z}-\hat{x}) \\
			\dot{\hat{z}}=A(t,t_{p})(\mathbb{L}\hat{y}-\mathbb{L}\hat{z}-\mathbb{L}\hat{u})\\
			\dot{\hat{u}}=A(t,t_{p})(\mathbb{L}\hat{y}-\mathbb{L}\hat{x}).
		\end{cases}
	\end{equation}	
	Thus, the origin of the system (\ref{AAb1})
	is equivalent to the equilibrium point of (\ref{AA1}). Therefore, we are going to analyze the convergence of the system (\ref{AAb1}).
	Motivated by \cite{YL-2023}, we apply the orthogonal transformation matrix:
	\begin{equation}\label{cedw}
		\begin{split}
			&\xi=[R_{1},R_{2}]^{T}\otimes{\rm I_{N}}\hat{x}={\rm col}[\xi_{1},\xi_{2}],\;\;\;
			\eta=[R_{1},R_{2}]^{T}\otimes{\rm I_{N}}\hat{y}={\rm col}[\eta_{1},\eta_{2}],\\
			&\delta=[R_{1},R_{2}]^{T}\otimes{\rm I_{N}}\hat{z}={\rm col}[\delta_{1},\delta_{2}],\;\;\;
			\omega=[R_{1},R_{2}]^{T}\otimes{\rm I_{N}}\hat{u}={\rm col}[\omega_{1},\omega_{2}],\\
		\end{split}
	\end{equation}
	where ${\rm I}_{N}$ denotes an $N$-dimensional identity matrix, $\xi_{1}$, $\eta_{1}$, $\delta_{1}$, $\omega_{1}\in \mathbb{R}$, $\xi_{2}$, $\eta_{2}$, $\delta_{2}$, $\omega_{2}\in \mathbb{R}^{N-1}$, $[R_{1},R_{2}]$ is an orthogonal matrix; $R_{1}=(1/\sqrt{N}) 1_{N}$ ( $1_{N}$ denotes an $N$-dimensional vector with the components being 1), $R_{2}^{T}R_{2}={\rm I}_{N-1}$, $R_{2}R_{2}^{T}={\rm I}_{N}-(1/N)1_{N}^{T}1_{N}$, $R_{1}^{T}R_{2}=0_{N-1}^{T}$. Thus, system (\ref{AAb1}) is
	rewritten as follows:
	\begin{equation*}
		\begin{cases}
			\dot{\xi}_{1}=A(t,t_{p})(-R_{1}^{T}h+\eta_{1}-\xi_{1}) \\
			\dot{\eta}_{1}=-A(t,t_{p})\xi_{1} \\
			\dot{\delta}_{1}=0\\
			\dot{\omega}_{1}=0,
		\end{cases}
		\begin{cases}
			\dot{\xi}_{2}=A(t,t_{p})(-R_{2}^{T}h+\eta_{2}-\xi_{2}) \\
			\dot{\eta}_{2}=A(t,t_{p})((R_{2}^{T}\mathbb{L}R_{2})(\frac{\varrho}{A(t,t_{p})}\delta_{2}-\eta_{2})-\xi_{2}) \\
			\dot{\delta}_{2}=A(t,t_{p})(R_{2}^{T}\mathbb{L}R_{2})(\eta_{2}-\delta_{2}-\omega_{2})\\
			\dot{\omega}_{2}=A(t,t_{p})(R_{2}^{T}\mathbb{L}R_{2})(\eta_{2}-\xi_{2}),
		\end{cases}
	\end{equation*}
	where $h=-\frac{\varrho}{A(t,t_{p})}(\nabla f(x)-\nabla f(x^{*}))$. 
	
	Now, our purpose is to prove that $\xi$ tends to the neighborhood of origin at the predefined time. Consider the following Lyapunov function
	\begin{equation}\label{L1}
		\begin{split}
			V(t)=&\frac{\rho}{2}(\xi^{T}\xi+\eta^{T}\eta)+\frac{\beta}{2}(\delta_{2}+\omega_{2})^{T}(\delta_{2}+\omega_{2})
			+\frac{\gamma}{2}(\xi-\eta)^{T}(\xi-\eta),
		\end{split}
	\end{equation}
	where $\rho$, $\beta$, and $\gamma$ are some positive constants.
	
	Note that 
	\begin{equation*}
		\begin{split}
			\xi^{T}\dot{\xi}
			=&A(t,t_{p})(-\hat{x}h-\xi^{T}\xi+\xi^{T}\eta),
			\\
			\eta^{T}\dot{\eta}=&
			A(t,t_{p})(-\eta^{T}\xi-\eta_{2}^{T}(R_{2}^{T}\mathbb{L}R_{2})\eta_{2})
			+\varrho(\eta_{2}^{T}(R_{2}^{T}\mathbb{L}R_{2})\delta_{2}),
			\\
			(\delta_{2}+\omega_{2})^{T}(\dot{\delta}_{2}+\dot{\omega}_{2})
			=&A(t,t_{p})( \delta_{2}+\omega_{2})^{T}
			(R_{2}^{T}\mathbb{L}R_{2})(2\eta_{2}-\delta_{2}-\omega_{2}
			-\xi_{2}),\\
			(\xi-\eta)^{T}(\dot{\xi}-\dot{\eta})
			=&A(t,t_{p}) (-\hat{x}h-\eta_{2}^{T}(R_{2}^{T}\mathbb{L}R_{2})\eta_{2}
			+\xi_{2}^{T}(R_{2}^{T}\mathbb{L}R_{2})\eta_{2}+\eta^{T}[R_{1},R_{2}]^{T}h\\
			&
			-\eta^{T}\eta+\eta^{T}\xi)
			-\varrho\eta_{2}^{T}(R_{2}^{T}\mathbb{L}R_{2})\delta_{2}
			-\varrho\xi_{2}^{T}(R_{2}^{T}\mathbb{L}R_{2})\delta_{2} .
		\end{split}
	\end{equation*}	
	
	Taking the time derivative of (\ref{L1}), we get that 
	\begin{equation}\label{de}
		\begin{split}
			\dot{V}(t)=&\rho(\xi^{T}\dot{\xi}+\eta^{T}\dot{\eta})+\beta(\delta_{2}+\omega_{2})^{T}(\dot{\delta}_{2}+\dot{\omega}_{2})
			+\gamma(\xi-\eta)^{T}(\dot{\xi}-\dot{\eta})\\
			=&\rho A(t,t_{p})(-\hat{x}^{T}h-\xi^{T}\xi-\eta_{2}^{T}(R_{2}^{T}\mathbb{L}R_{2})\eta_{2})
			-\rho\varrho\eta_{2}^{T}(R_{2}^{T}\mathbb{L}R_{2})\delta_{2}\\
			&+\beta A(t,t_{p})( \delta_{2}+\omega_{2})^{T}
			(-(R_{2}^{T}\mathbb{L}R_{2})(\delta_{2}+\omega_{2})
			+2(R_{2}^{T}\mathbb{L}R_{2})\eta_{2}\\
			&-(R_{2}^{T}\mathbb{L}R_{2})\xi_{2})
			+\gamma A(t,t_{p}) (-\hat{x}h
			-\eta_{2}^{T}(R_{2}^{T}\mathbb{L}R_{2})\eta_{2}
			+\xi_{2}^{T}(R_{2}^{T}\mathbb{L}R_{2})\eta_{2}\\
			&+\eta^{T}[R_{1},R_{2}]^{T}h
			-\eta^{T}\eta+\eta^{T}\xi)
			-\gamma\varrho\eta_{2}^{T}(R_{2}^{T}\mathbb{L}R_{2})\delta_{2}
			-\gamma\varrho\xi_{2}^{T}(R_{2}^{T}\mathbb{L}R_{2})\delta_{2}.	\end{split}
	\end{equation}
	Furthermore, according to Young’s inequality $x^{T}y\le (c/2)x^{T}y+ (1/2c)x^{T}y$, we obtain that
	\begin{equation*}
		\begin{split}
			&-A(t,t_{p})\hat{x}^{T}h
			\le\frac{\varrho(1+M^{2})}{2}\xi^{T}\xi,\;\;\;
			\eta_{2}^{T}(R_{2}^{T}\mathbb{L}R_{2})\eta_{2}
			\ge\lambda_{2}(\mathbb{L})\eta_{2}^{T}\eta_{2},\\
			&-\eta_{2}^{T}(R_{2}^{T}\mathbb{L}R_{2})\delta_{2}
			\le\frac{\lambda_{N}(\mathbb{L}^{2})}{2}\eta_{2}^{T}\eta_{2}+\frac{1}{2}\delta_{2}^{T}\delta_{2},\\
			&(\delta_{2}+\omega_{2})^{T}(R_{2}^{T}\mathbb{L}R_{2})(\delta_{2}+\omega_{2})
			\ge\lambda_{2}(\mathbb{L})(\delta_{2}+\omega_{2})^{T}(\delta_{2}+\omega_{2}),\\
			&(\delta_{2}+\omega_{2})^{T}(R_{2}^{T}\mathbb{L}R_{2})\xi_{2}
			\le\frac{\lambda_{N}(\mathbb{L}^{2})}{2\beta}(\delta_{2}+\omega_{2})^{T}(\delta_{2}+\omega_{2})+\frac{\beta}{2}\xi_{2}^{T}\xi_{2},\\
			&(\delta_{2}+\omega_{2})^{T}(R_{2}^{T}\mathbb{L}R_{2})\eta_{2}
			\le\frac{\lambda_{N}(\mathbb{L}^{2})}{2\beta}(\delta_{2}+\omega_{2})^{T}(\delta_{2}+\omega_{2})+\frac{\beta}{2}\eta_{2}^{T}\eta_{2},\\
			&\xi_{2}^{T}(R_{2}^{T}\mathbb{L}R_{2})\eta_{2}
			\le\frac{\gamma\lambda_{N}(\mathbb{L}^{2})}{2}\xi_{2}^{T}\xi_{2}+\frac{1}{2\gamma}\eta_{2}^{T}\eta_{2},\;\;\;
			\xi_{2}^{T}(R_{2}^{T}\mathbb{L}R_{2})\delta_{2}
			\le\frac{\lambda_{N}(\mathbb{L}^{2})}{2\gamma}\xi_{2}^{T}\xi_{2}+\frac{\gamma}{2}\delta_{2}^{T}\delta_{2},\\
			&\frac{1}{A(t,t_{p})}\eta^{T}[R_{1},R_{2}]^{T}h\le\frac{1}{2}\eta^{T}\eta+\frac{M^{2}}{2}\xi^{T}\xi^{T},\;\;\;
			\eta^{T}\xi\le\frac{1}{2\gamma}\eta^{T}\eta+\frac{\gamma}{2}\xi^{T}\xi,
		\end{split}
	\end{equation*}
	where $M=\max\{M_{i}\}$. Hence, we have that
	\begin{equation*}
		\begin{split}
			\dot{V}(t)
			\le&-A(t,t_{p})(\rho-\frac{\beta^{2}}{2}-\frac{\gamma^{2}\lambda_{N}(\mathbb{L}^{2})}{2}-\frac{\gamma^{2}}{2})\xi^{T}\xi
			-A(t,t_{p})(\rho\lambda_{2}(\mathbb{L})
			+\gamma\lambda_{2}(\mathbb{L})-\frac{1}{2}\\
			&-\beta^{2})\eta_{2}^{T}\eta_{2}
			-A(t,t_{p})(\beta\lambda_{2}(\mathbb{L})-\frac{3}{2}\lambda_{N}(\mathbb{L}^{2}))
			( \delta_{2}+\omega_{2})^{T}(\delta_{2}+\omega_{2})\\
			&-A(t,t_{p})(\gamma-\frac{1}{2})\eta^{T}\eta
			+\varrho\bigg(\frac{\gamma}{2}\eta^{T}\eta
			+\frac{(\rho+\gamma)(1+M^{2})+\gamma(\lambda_{N}(\mathbb{L}^{2})+M^{2})}{2}\xi^{T}\xi\\
			&+\frac{(\rho+\gamma)\lambda_{N}(\mathbb{L}^{2})}{2}\eta_{2}^{T}\eta_{2}+\frac{\rho+2\gamma}{2}\delta_{2}^{T}\delta_{2}\bigg).
		\end{split}
	\end{equation*}	
	Choose $\gamma\ge 2$, $\beta\ge(1+\lambda_{N}(\mathbb{L}^{2}))/(2\lambda_{2}(\mathbb{L}))$ and $\rho\ge\{(1+\beta+\gamma^{2}\lambda_{N}(\mathbb{L}^{2})+\gamma^{2})/2,(1-2\gamma\lambda_{2}(\mathbb{L}))/2\lambda_{2}(\mathbb{L})\}$ such that
	\begin{equation}\label{g1}
		\begin{split}
			&\frac{\gamma}{2}-\frac{1}{2}\ge\frac{1}{2},\;\;\;\;\;\;\;\;\;\;\;\;\;\;\;\;\;\;\;\;\;\;\;\;\;\;\;\;\;\;\;\;
			\beta\lambda_{2}(\mathbb{L})-\frac{3}{2}\lambda_{N}(\mathbb{L}^{2})\ge\frac{1}{2},\\
			&\rho-\frac{\beta}{2}-\frac{\gamma^{2}\lambda_{N}(\mathbb{L}^{2})}{2}-\frac{\gamma^{2}}{2}\ge\frac{1}{2},\;\;\;
			\rho\lambda_{2}(\mathbb{L})+\gamma\lambda_{2}(\mathbb{L})-\frac{1}{2}-\beta^{2}\ge0.
		\end{split}
	\end{equation}
	Note that
	\begin{equation}\label{lpjs}
		\begin{split}
			\frac{\min\{\rho,\beta\}}{2}{\rm col}[\xi,\eta,\omega_{2}]^{T}{\rm col}[\xi,\eta,\omega_{2}]		\le V(t)
			\le\frac{\rho+\beta+3\gamma}{2}{\rm col}[\xi,\eta,\delta_{2}+\omega_{2}]^{T}{\rm col}[\xi,\eta,\delta_{2}+\omega_{2}].
		\end{split}
	\end{equation}	
	It follows from (\ref{g1})-(\ref{lpjs}) that
	\begin{equation*}
		\begin{split}
			\dot{V}(t)\le&-\frac{1}{2}A(t,t_{p}){\rm col}[\xi,\eta,\delta_{2}+\omega_{2}]^{T}{\rm col}[\xi,\eta,\delta_{2}+\omega_{2}]
			+\varrho\Lambda_{1} {\rm col}[\xi,\eta,\omega_{2}]^{T}{\rm col}[\xi,\eta,\omega_{2}]\\
			\le&-\frac{A(t,t_{p})}{\rho_{1}+\beta+3\gamma}V(t)+\frac{\varrho\Lambda_{1}}{\min\{\rho,\beta\}}V(t),
		\end{split}
	\end{equation*}	
	where $\Lambda_{1}=\max\{[(\rho+\gamma)(1+M^{2})+\gamma(\lambda_{N}(\mathbb{L}^{2})+M^{2})]/2,[(\rho+\gamma)\lambda_{N}(\mathbb{L}^{2})+\gamma]/2,[\rho+2\gamma]/2\}$. 
	Set $\delta=\varrho\Lambda_{1}/\min\{\rho,\beta\}$, note that $\varrho$ is sufficient small such that $\varrho\le\alpha\min\{\rho,\beta\}/(D\Lambda_{1})$. By
	$
	\xi^{T}\xi=([R_{1},R_{2}]\otimes{\rm I}_{N}\cdot\hat{x})^{T}[R_{1},R_{2}]\otimes{\rm I}_{N}\cdot\hat{x}=\hat{x}^{T}\hat{x},
	$
	it follows from Theorem \ref{thr} that  
	$
	\parallel\Gamma(t,s)\parallel\le \mu_{\alpha-\delta D,D}(t,s).
	$
	This implies that
	$
	V(t)\le \mu_{\alpha-\delta D,D}(t,0)V(0).
	$
	By (\ref{L1}), it follows from Theorem \ref{vcx} that  
	$
	\parallel x(t)-x^{*}\parallel^{2}=\hat{x}^{T}\hat{x}=\xi^{T}\xi\le \mu_{\alpha-\delta D,D}(t,0)\frac{V(x(0))}{\rho/2}.
	$
	Hence,
	\begin{equation}\label{xishu2}
		\lim_{t\to t_{p+}}\parallel x(t)-x^{*}\parallel^{2}\le \mu_{\alpha-\delta D,2D/\rho}(t_{p},0)V(x(0)).
	\end{equation}
	Now let $\epsilon^{2}=\mu_{\alpha-\delta D,2D/\alpha}(t_{p},0)$, then
	$
	\lim_{t\to t_{p+}}\parallel x(t)-x^{*}\parallel^{2}\le \epsilon
	$.
	Similarly, we obtain that
	$	\lim_{t\to \infty}\parallel x(t)-x^{*}\parallel^{2}\le\epsilon$ for $\forall t\ge t_{p}$,
	and
	$
	\lim_{t\to \infty}\parallel x(t)-x^{*}\parallel^{2}\le0.
	$
	This complete the proof.
	
	Now we introduce the  Karush–Kuhn–Tucker (KKT) conditions. First, we first introduce the following lemmas. 
	\begin{lem}\label{kkt3}(\cite{GZ-2020})
		$x_{i}^{*}$ $(i\in P_{N})$ are the KKT points of problem (\ref{xd}) if and only if
		\begin{equation}\label{c1}
			\begin{cases}
				\nabla f_{i}(x_{i}^{*})=\nabla f_{j}(x_{j}^{*})\\
				\sum_{i=1}^{N}x^{*}_{i}=\sum_{i=1}^{N}q_{i}.
			\end{cases}
		\end{equation}
	\end{lem}
	
	Then the following corollaries are obtained
	Now we consider the cases that $f_{i}(x_{i})$ are strong convex, strictly convex, convex and non-convex.
	\begin{cor}\label{c14}
		Suppose that Assumption \ref{Ass1} holds, further assume that
		$f_{i}(x_{i})$ if continuously differentiable and $M$-smooth. Any of the following statement is correct:
		
		(1)	 MAS (\ref{A1}) converges to the KKT points of problem (\ref{xd}) at $t_{p}$.
		
		(2) If
		$f_{i}(x_{i})$ is strongly convex, then MAS  (\ref{A1}) converges to the unique optimal solution of problem (\ref{xd}) at $t_{p}$.
		
		(3) If
		$f_{i}(x_{i})$ ($\forall i\in P_{N}$) is strictly convex, then MAS  (\ref{A1}) converges to the unique optimal solution of problem (\ref{xd}) at $t_{p}$.
		
		(4) If $f_{i}(x_{i})$ ($\forall i\in P_{N}$) is convex, then MAS (\ref{A1}) converges to an optimal solution of problem (\ref{xd}) at $t_{p}$.
		
		(5)  MAS (\ref{A1}) converges to an local optimal solution of problem (\ref{xd}) at $t_{p}$.
	\end{cor}
	
	Now we consider the consensus-based distributed optimization problem.
	\begin{thm}\label{thxx}
		Suppose that Assumption \ref{Ass1} holds. Further assume that
		$f_{i}(x_{i})$ are continuously differentiable and $M$-smooth, then MAS (\ref{A2}) converges to its equilibrium points at $t_{p}$.
	\end{thm}
	
	{\it Proof:} 
	%Let $\nabla f(x)=[\nabla f_{1}(x_{1})^{T},f_{2}(x_{2})^{T},...,f_{N}(x_{N})^{T}]^{T}$, $\nabla f(x^{*})=[\nabla f_{1}(x_{1}^{*})^{T},f_{2}(x_{2}^{*})^{T},...,f_{N}(x_{N}^{*})^{T}]^{T}$, $x=[ x_{1}^{T},x_{2}^{T},...,x_{N}^{T}]^{T}$, $w=[ w_{1}^{T},w_{2}^{T},...,w_{N}^{T}]^{T}$, $w^{*}=[ w_{1}^{*T},w_{2}^{*T},...w_{N}^{*T}]^{T}$. Then the
	(\ref{A2}) can be rewritten in a compact form as follows
	\begin{equation}\
		\begin{cases}
			\dot{x}=A(t,t_{p})(-\frac{\varrho}{A(t,t_{p})}\nabla f(x)-\mathbb{L}x-\mathbb{L}w) \\
			\dot{w}=A(t,t_{p}) (\mathbb{L}w).
		\end{cases}
	\end{equation}
	Similarly, denote the equilibrium point of (\ref{AA1}) as $(x^{*}, w^{*})$ and
	let $\hat{x}=x-x^{*}$ and $\hat{w}=w-w^{*}$, then we have that
	\begin{equation*}
		\begin{cases}
			\dot{\hat{x}}=A(t,t_{p})(-\frac{\varrho}{A(t,t_{p})}(\nabla f(x)-\nabla f(x^{*}))-\mathbb{L}\hat{x}-\mathbb{L}\hat{w}) \\
			\dot{\hat{w}}=A(t,t_{p})\mathbb{L}\hat{x}.
		\end{cases}
	\end{equation*}	
	Through (\ref{cedw}), we obtain that 
	\begin{equation*}
		\begin{cases}
			\dot{\xi}_{1}=A(t,t_{p})(-R_{1}^{T}h-R_{1}^{T}\mathbb{L} \hat{x}) \\
			\dot{\eta}_{1}=0,
		\end{cases}
		\begin{cases}
			\dot{\xi}_{2}=A(t,t_{p})(-R_{2}^{T}h-R_{2}^{T}\mathbb{L}\hat{x}-(R_{2}^{T}\mathbb{L}R_{2})\eta_{2}) \\
			\dot{\eta}_{2}=A(t,t_{p})(R_{2}^{T}\mathbb{L}R_{2})\xi_{2},
		\end{cases}
	\end{equation*}
	where $h=-\frac{\varrho}{A(t,t_{p})}(\nabla f(x)-\nabla f(x^{*}))$.
	Consider the following Lyapunov function
	\begin{equation}\label{L2}
		\begin{split}
			V(t)=&\frac{\rho}{2}\xi^{T}\xi+\frac{\rho+\beta}{2}\eta_{2}^{T}\eta_{2}+\frac{\beta}{2}(\xi_{2}+\eta_{2})^{T}(\xi_{2}+\eta_{2}),
		\end{split}
	\end{equation}
	where $\rho$ and $\beta$ are positive constants.
	Taking the time derivative of (\ref{L2}), we obtain that 
	\begin{equation}\label{v2d}
		\begin{split}
			\dot{V}(t)
			=&\rho(\xi^{T}\dot{\xi}+\eta_{2}^{T}\dot{\eta}_{2})+(\rho+\beta)\eta_{2}^{T}\dot{\eta}_{2}+\beta(	\xi_{2}^{T}+\eta_{2}^{T})(\dot{\xi}_{2}^{T}+\dot{\eta}_{2}^{T})\\
			=&\rho A(t,t_{p})(-\hat{x}^{T}\mathbb{L}\hat{x})
			+\beta A(t,t_{p})(-\eta_{2}^{T}(R_{2}^{T}\mathbb{L}R_{2})\eta_{2})\\
			&-\beta\varrho(\xi_{2}^{T}+\eta_{2}^{T})R_{2}^{T}(\nabla f(x)-\nabla f(x^{*}))
			-\rho\varrho\hat{x}(\nabla f(x)-\nabla f(x^{*}))
		\end{split}
	\end{equation}
	By Young’s inequality, we obtain that 
	\begin{equation}\label{ggg2}
		\begin{split}
			&-\hat{x}^{T}\mathbb{L}\hat{x}\ge -\lambda_{2}(N)\hat{x}^{T}\hat{x}=-\lambda_{2}(\mathbb{L})\xi^{T}\xi,\\
			&-\eta_{2}^{T}(R_{2}^{T}\mathbb{L}R_{2})\eta_{2}\le-\lambda_{2}(\mathbb{L})\eta_{2}^{T}\eta_{2},\\
			&-\hat{x}(\nabla f(x)-\nabla f(x^{*}))\le\frac{1+M^{2}}{2}\xi^{T}\xi,\\
			&	-(\xi_{2}^{T}+\eta_{2}^{T})R_{2}^{T}(\nabla f(x)-\nabla f(x^{*}))
			\le\frac{M^{2}+\lambda_{2}(\mathbb{L})}{2}\xi^{T}\xi+\frac{\lambda_{2}(\mathbb{L})}{2}\eta_{2}^{T}\eta_{2}.
		\end{split}
	\end{equation}
	Therefore, we have that 
	\begin{equation*}
		\begin{split}
			\dot{V}(t)\le&-A(t,t_{p})\rho\lambda_{2}(\mathbb{L})\xi^{T}\xi-A(t,t_{p})\beta\lambda_{2}(\mathbb{L})\eta_{2}^{T}\eta_{2}\\
			&+\varrho\bigg(\frac{\rho(1+M^{2})+\beta(M^{2}+\lambda_{2}(\mathbb{L}))}{2}\xi^{T}\xi
			+\frac{\beta\lambda_{2}(\mathbb{L})}{2}\eta_{2}^{T}\eta_{2}\bigg).
		\end{split}
	\end{equation*}
	Choose $\rho\ge 1/2\lambda_{2}(\mathbb{L})$ and $\beta\ge1/2\lambda_{2}(\mathbb{L})$ such that
	\begin{equation}\label{gg1}
		\begin{split}
			\rho\lambda_{2}(\mathbb{L})\ge\frac{1}{2},\;\;\;
			\beta\lambda_{2}(\mathbb{L})\ge\frac{1}{2}.
		\end{split}
	\end{equation}	
	Then 
	$
	\dot{V}(t)\le-\frac{1}{2}A(t,t_{p})(\xi^{T}\xi+\eta_{2}^{T}\eta_{2})
	+\varrho\widetilde{\Upsilon}_{1}(\xi^{T}\xi+\eta_{2}^{T}\eta_{2}),
	$
	where $\Upsilon_{1}=\max\{[\rho(1+M^{2})+\beta(M^{2}+\lambda_{2}(\mathbb{L}))]/2, \beta\lambda_{2}(\mathbb{L})/2\}$.
	Note that
	\begin{equation*}
		\begin{split}
			\frac{\rho+3\beta}{2}{\rm col}[\xi,\eta_{2}]^{T}{\rm col}[\xi,
			\eta_{2}]
			\ge V(t)
			\ge\frac{\min\{\rho,\beta\}}{2}{\rm col}[\xi,\eta_{2}]^{T}{\rm col}[\xi,
			\eta_{2}].
		\end{split}
	\end{equation*}
	Hence, we obtain that
	$
	\dot{V}(t)\le-\frac{A(t,t_{p})}{\rho+3\beta}(\xi^{T}\xi+\eta_{2}^{T}\eta_{2})
	+\frac{\varrho\Upsilon_{1}}{\min\{\rho,\beta\}}(\xi^{T}\xi+\eta_{2}^{T}\eta_{2}).
	$
	Similar to the proof of Theorem (\ref{thth1}), we obtain the desired results. 
	
	Here we introduce a important lemma similar to Lemma \ref{kkt3}: 
	\begin{lem}(\cite{YL-2023})\label{cc2}
		$x_{i}^{*}$ $(i\in P_{N})$ are the KKT points of problem (\ref{xx}) if and only if
		\begin{equation}\label{c2}
			\begin{cases}
				0=-\frac{\varrho}{A(t,t_{p})}\nabla f(x^{*})-\mathbb{L}x^{*}-\mathbb{L}w^{*} \\
				0=\mathbb{L}x^{*}.
			\end{cases}
		\end{equation}
	\end{lem}
	
	Similar to Corollary \ref{c14}, according to Lemma \ref{cc2} we have that
	
	\begin{cor}\label{c9}
		Suppose that Assumption \ref{Ass1} holds, further assume that
		$f_{i}(x_{i})$ if continuously differentiable and $M$-smooth. Any of the following statement is correct:
		
		(1) MAS (\ref{A2}) converges to the KKT points of problem (\ref{xx})  at $t_{p}$.
		
		(2) If
		$f_{i}(x_{i})$ is strongly convex, then MAS (\ref{A2}) converges to the unique optimal solution of problem (\ref{xx}) at $t_{p}$.
		
		(3) If
		$f_{i}(x_{i})$ is strictly convex, then MAS (\ref{A2}) converges to the unique optimal solution of problem (\ref{xx}) at  $t_{p}$.
		
		(4) If
		$f_{i}(x_{i})$ is convex, then MAS (\ref{A2}) converges to an optimal solution of problem (\ref{xx}) at $t_{p}$.
		
		(5) MAS (\ref{A2}) converges to an optimal solution of problem (\ref{xx}) at $t_{p}$.
	\end{cor}
	Now we consider the generalized smoothness.
	\begin{thm}\label{ehn1}
		Suppose that Assumption \ref{Ass1} holds. Further assume that
		$f_{i}(x_{i})$ is continuously differentiable and generalized smooth, then  MAS (\ref{A1}) converges to its equilibrium points at $t_{p}$.
	\end{thm}
	
	{\it Proof:} The previous proof is similar to Theorem \ref{thth1}, we directly start calculating the Lyapunov function.
	Note that
	\begin{equation}\label{xhmm}
		\begin{split}
			-A(t,t_{p})\hat{x}^{T}h
			=&-\varrho\hat{x}(\nabla f(x)-\nabla f(x^{*}))
			\le\varrho\bigg(\frac{M^{2}+2}{2}\xi^{T}\xi+\frac{\widetilde{M}^{2}M^{2}+\widetilde{M}^{2}}{2}\bigg),
		\end{split}				
	\end{equation}
	where $M=\max\{M_{i}\}$ and $\widetilde{M}=\max\{\widetilde{M}_{i}\}$. We have
	$
	\frac{1}{A(t,t_{p})}\eta^{T}[R_{1},R_{2}]^{T}h
	\le\frac{1}{2}\eta^{T}\eta+\frac{M^{2}+1}{2}\xi^{T}\xi+\frac{\widetilde{M}^{2}M^{2}+\widetilde{M}^{2}}{2},
	$
	Then (\ref{de}) can be written as 
	\begin{equation*}
		\begin{split}
			\dot{V}(t)
			\le&-A(t,t_{p})(\rho-\frac{\beta^{2}}{2}-\frac{\gamma^{2}\lambda_{N}(\mathbb{L}^{2})}{2}-\frac{\gamma^{2}}{2})\xi^{T}\xi
			-A(t,t_{p})(\rho\lambda_{2}(\mathbb{L})+\gamma\lambda_{2}(L)\\
			&-\frac{1}{2}-\beta^{2})\eta_{2}^{T}\eta_{2}
			-A(t,t_{p})(\beta\lambda_{2}(L)-\frac{3\lambda_{N}(\mathbb{L}^{2})}{2})
			( \delta_{2}+\omega_{2})^{T}(\delta_{2}+\omega_{2})\\
			&-A(t,t_{p})(\gamma-\frac{1}{2})\eta^{T}\eta
			+\varrho\bigg(\frac{\gamma}{2}\eta^{T}\eta
			+\frac{\gamma\lambda_{N}(\mathbb{L}^{2})+(M^{2}+2)(\rho+\gamma)}{2}\xi^{T}\xi\\
			&+\frac{\gamma(M^{2}+1)\xi^{T}\xi}{2}+\frac{(\rho+\gamma)\lambda_{N}(\mathbb{L}^{2})}{2}\eta_{2}^{T}\eta_{2}+\frac{\rho+2\gamma}{2}\delta_{2}^{T}\delta_{2}\bigg)\\
			&+\varrho\bigg(\frac{(\rho+2\gamma)(\widetilde{M}^{2}M^{2}+\widetilde{M}^{2})}{2}\bigg).
		\end{split}
	\end{equation*}
	It follows from(\ref{g1})-(\ref{lpjs}) that
	$
	\dot{V}(t)
	\le-\frac{A(t,t_{p})}{\rho+\beta+3\gamma}V(t)+\frac{\varrho\Lambda_{2}}{\min\{\rho,\beta\}}V(t)+\varrho\Lambda_{3},
	$
	where $\Lambda_{2}=\max\{[\gamma\lambda_{N}(\mathbb{L}^{2})+(M^{2}+2)(\rho+\gamma)+\gamma(M^{2}+1)]/2,[(\rho+\gamma)\lambda_{N}(\mathbb{L}^{2})+\gamma]/2,[\rho+2\gamma]/2\}$ and $\Lambda_{3}=(\rho+2\gamma)(\widetilde{M}^{2}M^{2}+\widetilde{M}^{2})/2$. 
	
	Note that $\varrho$ is sufficient small such that 
	$
	\frac{\varrho\Lambda_{2}}{\min\{\rho,\beta\}}\le\alpha/ D.
	$
	Thus, according to Theorem \ref{thr},
	$
	\dot{V}(t)
	=-\frac{A(t,t_{p})}{\rho+\beta+3\gamma}V(t)+\frac{\varrho\Lambda_{2}}{\min\{\rho,\beta\}}V(t)
	$
	also admits a predefined-time $\mu$ contraction.
	Furthermore, by Theorem $\ref{AVFV}$, $V(t)$ is bounded and achieves predefined-time approximate convergence to the bounded neighborhood of $0$ at $t_{p}$. Note that generalized smoothness alters the convergence rate, but the system still converges to the optimal solution in essence. This completes the proof.
	
	\begin{thm}\label{tha22}
		Suppose that Assumption \ref{Ass1} holds. Further assume that
		$f_{i}(x_{i})$ is continuously differentiable and generalized smooth, then MAS (\ref{A2}) converges to its equilibrium points at $t_{p}$.
	\end{thm}
	
	{\it Proof:} The previous proof is similar to Theorem \ref{thxx}, we directly start calculating the Lyapunov function.
	Combining (\ref{v2d})-(\ref{ggg2}), (\ref{xhmm}) and
	\begin{equation*}
		\begin{split}
			&(\xi_{2}^{T}+\eta_{2}^{T})R^{T}(\nabla f(x)-\nabla f(x^{*}))\\
			\le&\frac{\lambda_{2}(\mathbb{L})}{2}(\xi_{2}^{T}\xi_{2}+\eta_{2}^{T}\eta_{2})+\frac{M^{2}+1}{2}\xi^{T}\xi+\frac{\widetilde{M}^{2}M^{2}+\widetilde{M}^{2}}{2},
		\end{split}	
	\end{equation*}
	we obtain that
	\begin{equation*}
		\begin{split}
			\dot{V}(t)\le&-A(t,t_{p})\rho\lambda_{2}(\mathbb{L})\xi^{T}\xi-A(t,t_{p})\beta\lambda_{2}(\mathbb{L})\eta_{2}^{T}\eta_{2}\\
			&+\varrho\bigg(\frac{\rho(2+M^{2})+\beta(M^{2}+1+\lambda_{2}(\mathbb{L}))}{2}\xi^{T}\xi\\
			&+\frac{\beta\lambda_{2}(\mathbb{L})}{2}\eta_{2}^{T}\eta_{2}\bigg)
			+\varrho\frac{(\rho+\beta)(\widetilde{M}^{2}M^{2}+\widetilde{M}^{2})}{2}.
		\end{split}
	\end{equation*}
	It follows from (\ref{gg1}) that
	$
	\dot{V}(t)
	\le-\frac{A(t,t_{p})}{\rho+\beta+3\gamma}V(t)+\frac{\varrho\Upsilon_{2}}{\min\{\rho,\beta\}}V(t)+\varrho\Upsilon_{3},
	$
	where $\Upsilon_{2}=\max\{[\rho(M^{2}+2)+\beta(M^{2}+1+\lambda_{2}(\mathbb{L}))]/2,\beta\lambda_{2}(\mathbb{L})\}$ and $\Upsilon_{3}=(\rho+\beta)(\widetilde{M}^{2}M^{2}+\widetilde{M}^{2})/2$. Similar to Theorem \ref{ehn1}, we obtain the desired results.

	Then we have the following corollary.
	\begin{cor}\label{C7} 
		Replace $M$-smooth by generalized smooth, the conclusions of Corollary \ref{c14}-Corollary \ref{c9} also holds.
	\end{cor}
	Strongly convex is a necessary for predefined-time optimization in previous works. In some examples, if cost functions are non-convex, the original MASs in \cite{YL-2023,GZ-2020} cannot converge to local optimal points, whereas the multi-agent system we propose precisely addresses this issue.
	\begin{rem}\label{re1234}
		(1) It is noted that although the MASs we designed are different from the original MASs in existing papers, such as \cite{YL-2023} and \cite{GZ-2020}, the TBGs are still applicable to the previous MASs. Namely, the original MASs can still achieve predefined-time optimization under strong convexity if using TBG in Definition \ref{aptu}. However, strongly convex is a necessary for predefined-time optimization. If cost functions are non-convex, the original MASs in \cite{YL-2023,GZ-2020} under new TBGs cannot converge to local optimal points, whereas the MASs we propose precisely addresses this issue.
		
		(2) Two proposed MASs achieve predefined-time optimization under strong convexity. This imply that our results can achieve the same effect as conclusions in \cite{YL-2023} and \cite{GZ-2020}. Different from the MASs provided in \cite{YL-2023,KZ-2023,PW-2017,GZ-2018,GZ-2020,LS-2022,HWW-2022,PWJ-2017,MA-2022,BQZA-2019,GC-2022},  Theorem \ref{thth1}-Theorem \ref{tha22} show that the proposed MASs in this paper converge to their equilibrium points or the neighborhood of equilibrium points at a predefined time in the case of whether the cost function is strongly convex or not. Additionally, combine the Particle Swarm Optimization method \cite{ZJJ-2022}, the proposed MASs can converge to the global optimal solution with a probability of 1. Detailed examples can be found in the Section \ref{s5}.
		
		(3) In certain specific scenarios, the utilization of  smoothness conditions can lead to faster convergence compared to using traditional smoothness conditions, primarily due to the reduction of Lipschitz coefficients. Such specific scenarios occur when functions are continuously differentiable in certain neighborhood but exhibit large Lipschitz constants ($M$ is sufficiently large and $\widetilde{M}$ is sufficiently small). Furthermore, note that the value of $x$ at $t_{p}$ does not exceed a specific bound of optimal solution, i.e., $\lim_{t\to+\infty}\parallel x(t)-x^{*}\parallel\le \epsilon$. Therefore, under the condition that some coefficients are allowed ($M$ is sufficiently large and $\widetilde{M}$ is sufficiently small), the MASs under generalized smoothness exhibit faster convergence rate and CPU speed than the MASs under smoothness. Detailed examples can be found in the Section \ref{s5}.
	\end{rem}
	\begin{rem}\label{re78}
		As an example, let's design a MAS in form of (\ref{A1}) to solve RAP.
		In conducting specific numerical computations, the primary step involves determining the value of $\varepsilon$, i.e. a neighborhood of the optimal solution that arrives within a predefined time. Then we select proper coefficients $\gamma$, $\beta$ and $\rho$ such that inequalities (\ref{g1}) holds. In next step, select proper coefficients $\varrho$, and $A(t,t_{p})$ such that Theorem \ref{thr} and inequalities (\ref{xishu2}) hold. In $\mathcal{B}2)$, we note that $A(t,t_{p})$ can be a sufficiently large constant (we also give a example that $A(t,t_{p})$ is a function). This is consistent with the $(17b)$ and $(17c)$ stated in \cite{YL-2023}. Furthermore, Example 1 in \cite{YL-2023} also validates this statement. Detailed examples can be found in the Section \ref{s5}.
	\end{rem}
	
	\section{Illustrated Examples}\label{s5}
	In this section, to prove the effectiveness of two MASs, the obtained theoretical results are tested by the a large number of examples. Note that the proposed TBGs in this paper can be also applied into the MASs from \cite{YL-2023} and \cite{GZ-2020}. Hence we utilize the TBG in Definition \ref{aptu} in all MASs. However, the proposed MASs in this paper are significantly different from the Previous MASs, not only in terms of the inclusion of TBG, but also in the substantial changes made to the differential equation system. Hence, in order to highlight this point, we choose the simplest form of TBGs: a constant, which similar to the TBG3 in Example 1 form \cite{YL-2023}. We give the form of TBG when it's less than $t_{p}$. Additionally, the selection of $A(t,t_{p})$ and $\varrho$ can be derived following the methodology outlined in Remark \ref{re78}. The optimization methods are simulated in MATLAB R2017b and run on an Intel(R) Core(TM) i5-9750H CPU at 2.60 GHz, with an NVIDIA GeForce GTX 1660 Ti 6 GB graphics card.
	\subsection{Examples of RAPs}
	First, we are going to show that MASs (\ref{A1}) can achieve comparable performance to the original MASs for RAP in \cite{GZ-2020} under the assumptions of (2) in Corollary \ref{c14}.
	\begin{exmp}\label{E21}
		This example consider the widely discussed economic dispatch problems. We study the cost functions with the form of derivatives:	$\nabla f_{1}(x_{1})=4x_{1}+3$, $\nabla f_{2}(x_{2})=2x_{2}+4$, $\nabla f_{3}(x_{3})=x_{3}+5$,
		and $\nabla f_{4}(x_{4})=3x_{4}+2$.
		We set
		$x_{1}(0)=40$, $x_{2}(0)=35$, $x_{3}(0) = 45$, $x_{4}(0) = 40$, 
		$y_{i}(0)=0$, $z_{i}(0)=0$,
		$u_{i}(0)=0$ ($i\in P_{4}$). And   $L=[2,-1,0,-1;-1,2,-1,0;0,-1,2,-1;-1,0,-1,2]$.
		Set the total demand of this system as $q_{0} = 145$, $t_{p}=7$.
		Now we design the TBG as $A(t,t_{p})=2$ ($t\le t_{p}$) in original MASs. Meanwhile we set $A(t,t_{p})=80$ ($t\le t_{p}$) and $\varrho=8$ in MAS (\ref{A1}). From Fig. \ref{f4} (a), the states converge to a state, which is extremely close to the optimal solution at $t_{p}=7$. Fig. \ref{f4} (a) shows that MAS (\ref{A1}) is also applicable under strong convexity conditions.
		\begin{figure}[htph]
			\centering
			\subfigure[``old" and ``new" represent the transient behaviors under original MAS and  MAS (\ref{A1}) in Example \ref{E21}, respectively. ]{\includegraphics[width=4cm,height=4cm]{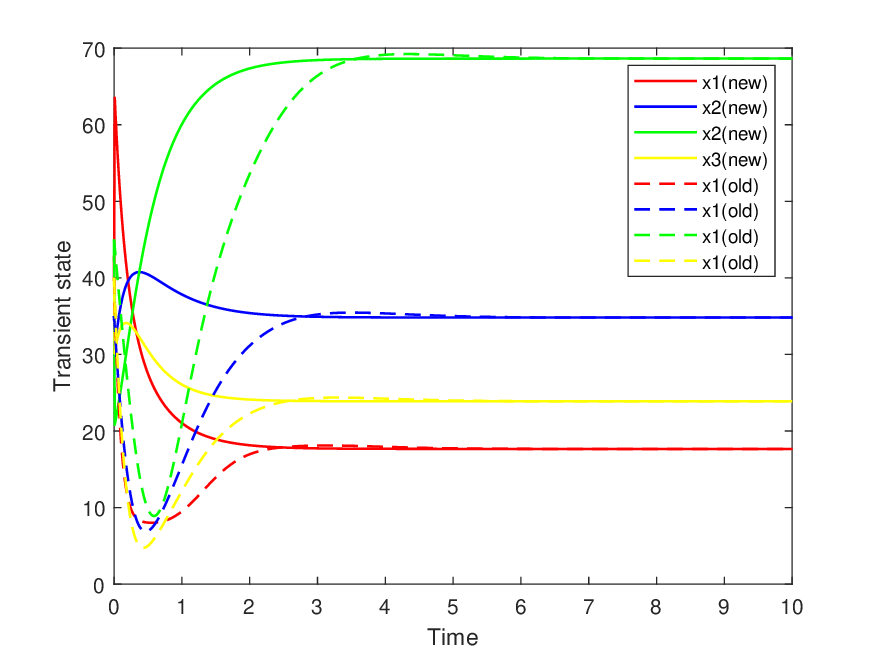}}
			\subfigure[Transient behaviors of $x_{i}(t)$ under original MAS in Example \ref{5e6}. ]{\includegraphics[width=4cm,height=4cm]{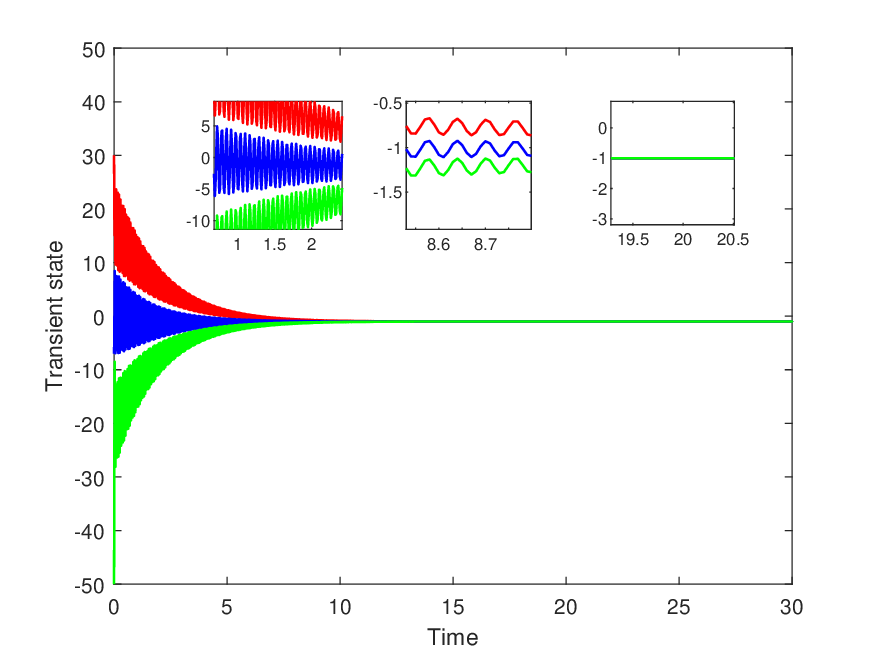}}
			\subfigure[Transient behaviors of $x_{i}(t)$ under MAS (\ref{A1}) in Example \ref{5e6}.]{\includegraphics[width=4cm,height=4cm]{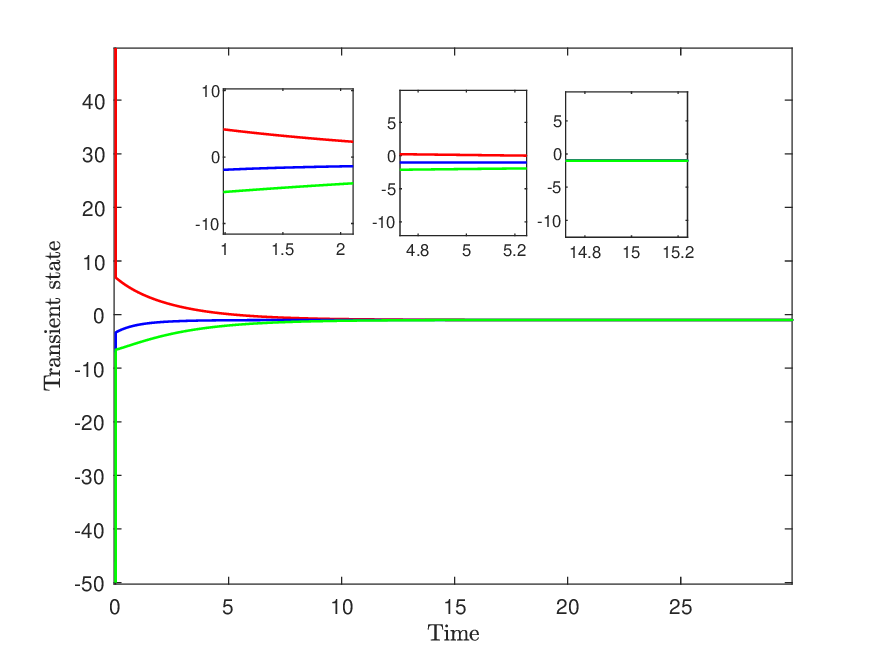}}	
			\caption{Transient behaviors of $x_{i}(t)$ in the different MASs.}\label{f4}
		\end{figure}
	\end{exmp}
	%%%%%%%%%%%%%%%%%%%%%%%%%%%%%%%%%%%%%%%%%%%%%%%%%5
	
	Now we aim to study the transient state of $x_{i}(t)$ for proposed MAS (\ref{A1}) and the original MAS under the assumptions of (3) in Corollary \ref{c14}.
	\begin{exmp}\label{5e6}
		Consider the cost functions with the form of derivatives: 
		$\nabla f_{i}(x_{i})=0.011x_{i}(x_{i}\le -1);$
		$\nabla f_{i}(x_{i})=(x_{i}+0.9)^{3}-0.01(-1\le x_{i}\le -0.8)$;
		$\nabla f_{i}(x_{i})=(x_{i}+0.7)^{3}-0.008, (-0.8\le x_{i}\le -0.6)$;
		$\nabla f_{i}(x_{i})=(x_{i}+0.5)^{3}-0.006,  (-0.6\le x_{i}\le -0.4)$;
		$\nabla f_{i}(x_{i})=(x_{i}+0.3)^{3}-0.004,  (-0.4\le x_{i}\le -0.2)$;
		$\nabla f_{i}(x_{i})=(x_{i}+0.1)^{3}-0.002,  (-0.2\le x_{i}\le 0)$;
		$\nabla f_{i}(x_{i})=(x_{i}-0.1)^{3}     ,   (0\le x_{i}\le 0.2)$;
		$\nabla f_{i}(x_{i})=(x_{i}-0.3)^{3}+0.002,  (0.2\le x_{i}\le 0.4)$;
		$\nabla f_{i}(x_{i})=(x_{i}-0.5)^{3}+0.004,  (0.4\le x_{i}\le 0.6)$;
		$\nabla f_{i}(x_{i})=(x_{i}-0.7)^{3}+0.006,  (0.6\le x_{i}\le 0.8)$;
		$\nabla f_{i}(x_{i})=(x_{i}-0.9)^{3}+0.008,  (0.8\le x_{i}\le 1$;
		$\nabla f_{i}(x_{i})=0.009x_{i}, (1\le x_{i})$.
		We set 
		$x_{1}(0)=50$, $x_{2}(0)=-7$, $x_{3}(0) = -50$,  
		$y_{i}(0)=0$,
		$z_{i}(0)=0$,
		$u_{i}(0)=0$ ($i\in P_{3}$),  $L=[1,-1,0;-1,2,-1;0,-1,1]$.
		We set $q_{0} = -3$ and $t_{p}=15$. 
		Now we design the TBG as $A(t,t_{p})=80$ ($t\le t_{p}$) in original MASs. Meanwhile we set $A(t,t_{p})=800$ ($t\le t_{p}$) and $\varrho=240$ in MAS (\ref{A1}).
		Note that we consider the same cost function for $i=1$, $2$, $3$, because it allowed us to observe the optimal solution in a more straightforward and simplified manner ($x^{*}$=[-1;-1;-1]). Fig. \ref{f4} (b) (original MAS) and Fig. \ref{f4} (c) (MAS (\ref{A1})) shows that the states converge to a state that is extremely close to the optimal solution at $t_{p}=15$. However, Fig. \ref{f4} (b) shows that the convergence process of $x$ exhibits an oscillatory pattern. And Fig. \ref{f4} (c) shows a relatively smooth trend.
	\end{exmp}
	
	Now we aim to show that under the assumptions stated of (4) in Corollary \ref{c14}, the transient state of $x_{i}(t)$ for proposed MAS (\ref{A1}) and the original MAS, respectively. 
	\begin{exmp}\label{5e7}
		Consider the cost functions with form of derivatives
		\begin{equation}
			\nabla
			f_{i}(x_{i})=\begin{cases}
				0.8\sin\frac{1}{0.4}-\cos\frac{1}{0.4}&x_{i}\ge0.4;\\
				2x\sin\frac{1}{x_{i}}-\cos\frac{1}{x_{i}}&0.2<x_{i}<0.4;\\
				0.4\sin\frac{1}{0.2}-\cos\frac{1}{0.2}&x_{i}\le0.2;
			\end{cases}
		\end{equation}
		Consider 
		$x_{1}(0)=5$, $x_{2}(0)=-5$, $x_{3}(0) = -1$,  
		$y_{i}(0)=0$ ($i\in P_{3}$),
		$z_{i}(0)=0$ ($i\in P_{3}$),
		$u_{i}(0)=0$ ($i\in P_{3}$),  $L=[1,-1,0;-1,2,-1;0,-1,1]$.
		Set $q_{0} = 3$ and $t_{p}=7$. 
		Now we design the TBG as $A(t,t_{p})=2000$ ($t\le t_{p}$) in original MAS. Meanwhile we set $A(t,t_{p})=8000$ ($t\le t_{p}$) and $\varrho=80$ in MAS (\ref{A1}).
		Fig. \ref{f5} (a) shows that the variation process of $x$ exhibits an oscillatory pattern and fails to converge. However, Fig. \ref{f5} (b) shows that the states converge to a state that is extremely close to the optimal solution at $t_{p}=7$ and the convergence process of $x$ exhibits a relatively smooth trend. Note that, due to the convexity of cost functions, the optimal solution is not unique. Hence, the MAS converges to one of the optimal solutions of the optimization problem.
		
		\begin{figure}[htph] 
			\centering
			\subfigure[Transient behaviors of $x_{i}(t)$ under original MAS in Example \ref{5e7}.  ]{\includegraphics[width=3cm,height=3cm]{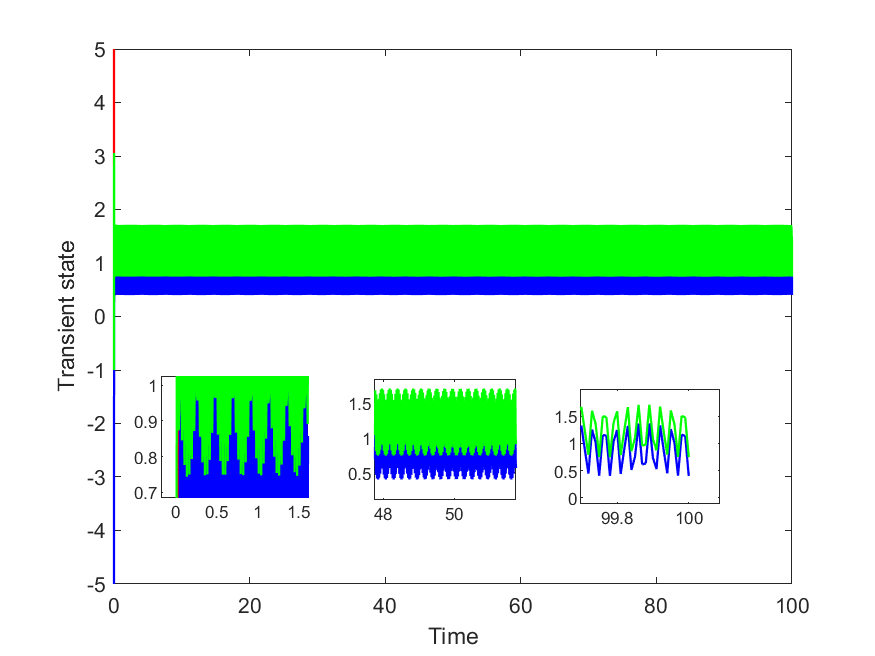}}
			\subfigure[Transient behaviors of $x_{i}(t)$ under MAS (\ref{A1}) in Example \ref{5e7}.]{\includegraphics[width=3cm,height=3cm]{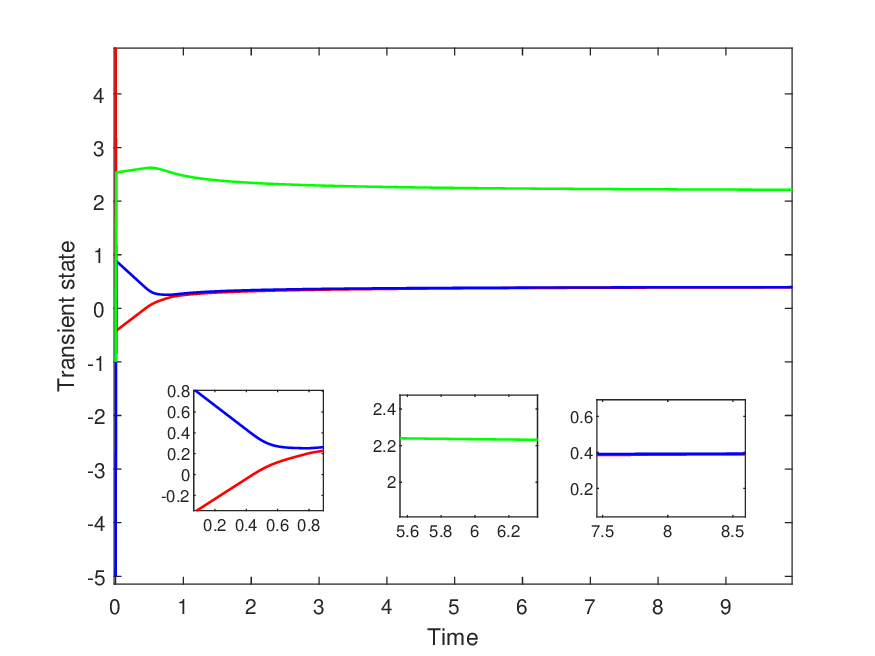}}	
			\subfigure[Transient behaviors of $x_{i}(t)$ under original MAS in Example \ref{5e8}. ]{\includegraphics[width=3cm,height=3cm]{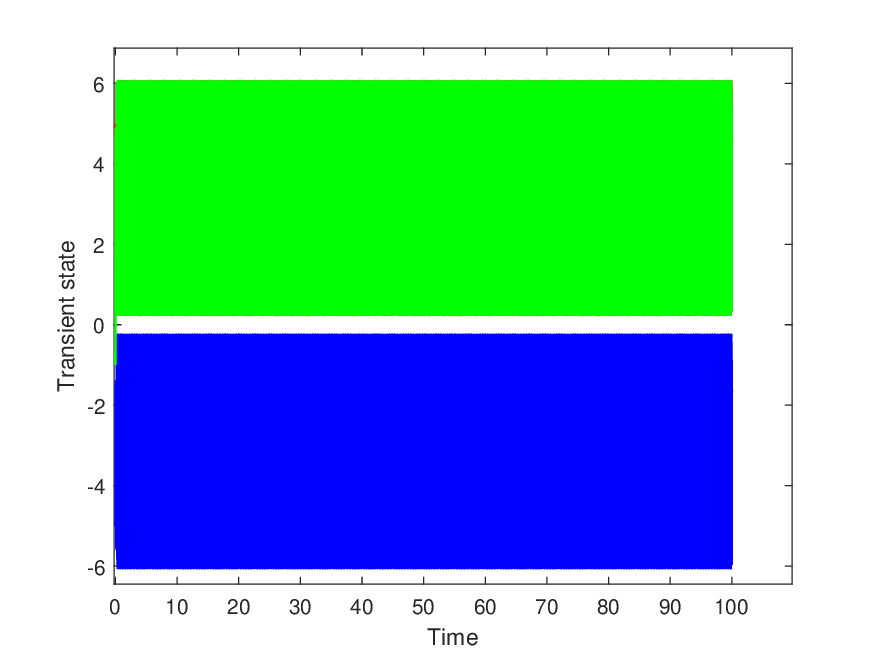}}
			\subfigure[Transient behaviors of $x_{i}(t)$ under MAS (\ref{A1}) in Example \ref{5e8}.]{\includegraphics[width=3cm,height=3cm]{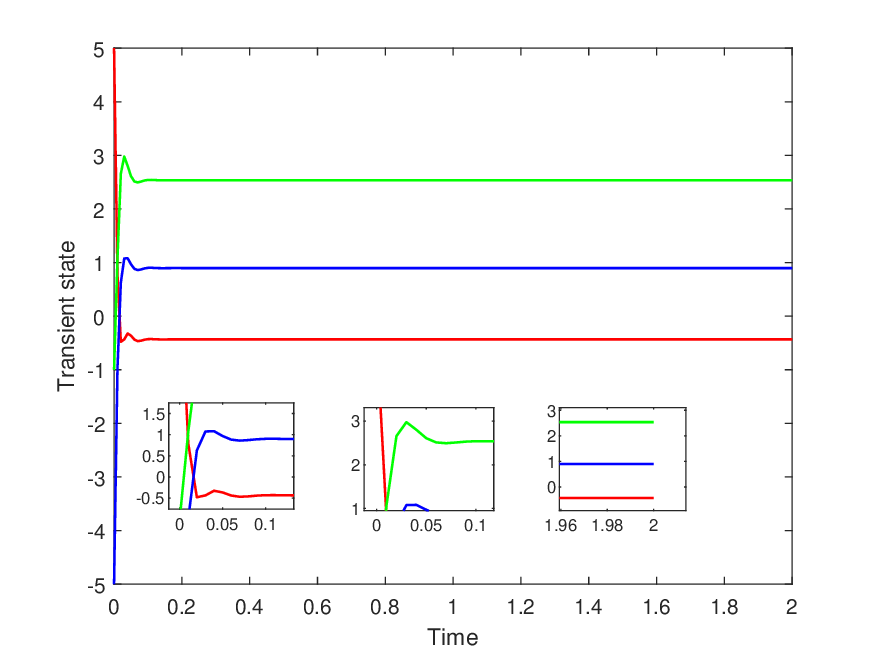}}
			\caption{Transient behaviors of $x_{i}(t)$ under the different MASs.}\label{f5}
		\end{figure}
		
	\end{exmp}
	
	Now we aim to show that under the assumptions of (5) in Corollary \ref{c14}, the transient state of $x_{i}(t)$ for proposed MAS (\ref{A1}) and the original MAS. 
	
	\begin{exmp}\label{5e8}
		We consider the cost functions:
		$
		f_{i}(x_{i})=x_{i}^{2}\sin\frac{1}{x_{i}}
		$
		Consider 
		$x_{1}(0)=5$, $x_{2}(0)=-5$, $x_{3}(0) = -1$,  
		$y_{i}(0)=0$ ($i\in P_{3}$), 
		$z_{i}(0)=0$ ($i\in P_{3}$), 
		$u_{i}(0)=0$ ($i\in P_{3}$),   $L=[1,-1,0;-1,2,-1;0,-1,1]$.
		Set $q_{0} = 3$ and $t_{p}=1.8$.
		Now we design the TBG as $A(t,t_{p})=200$ ($t\le t_{p}$) in original MASs. Meanwhile we set $A(t,t_{p})=100$ ($t\le t_{p}$) and $\varrho=1$ in MAS (\ref{A1}).
		Fig. \ref{f5} (c) shows that the variation process of $x$ exhibits an oscillatory pattern and fails to converge. However, Fig. \ref{f5} (d) shows that the states converge to a state that is extremely close to the optimal solution at $t_{p}=1.8$ and the convergence process of $x$ exhibits a relatively smooth trend. Note that, due to the convexity of cost functions, the optimal solution is not unique. Hence, the MAS converges to one of the optimal solutions of the optimization problem.

	\end{exmp}	
	
	By replacing $M$-smooth by  smooth in (2) of Corollary \ref{c14}, we aim to show the transient state of $x_{i}(t)$ for proposed MAS (\ref{A1}) and the original MAS.		
	
	\begin{exmp}	\label{5e9}	
		We consider the cost functions with form of derivatives:
		$\nabla	f_{1}(x_{1})=100x_{1}+90 ( x_{1}\ge\frac{1}{10})$;
		$\nabla	f_{1}(x_{1})=1000x_{1}    (-\frac{1}{10}<x_{1}<\frac{1}{10})$;
		$\nabla	f_{1}(x_{1})=100x_{1}-90 (x_{1}\le-\frac{1}{10})$.	
		$\nabla	f_{2}(x_{2})=2\nabla f_{1}(x_{2})$ and $\nabla	f_{3}(x_{3})=3\nabla f_{1}(x_{3})$. Easy to find that $|\nabla f_{1}(x_{1})-\nabla	f_{1}(y_{1})|\le 1000|x_{1}-y_{1}|$ and $|\nabla f_{1}(x_{1})-\nabla f_{1}(y_{1})|\le 100|x_{1}-y_{1}|+200$ for all $x_{1}$, $y_{1}\in\mathbb{R}$.
		Set 
		$x_{1}(0)=50$, $x_{2}(0)=-7$, $x_{3}(0) =- 50$,  
		$y_{i}(0)=0$ ($i\in P_{3}$),
		$z_{i}(0)=0$ ($i\in P_{3}$),
		$u_{i}(0)=0$ ($i\in P_{3}$),
		$L=[1,-1,0;-1,2,-1;0,-1,1]$.
		Set $q_{0} = 3$ and $t_{p}=70$. 
		Now we design the TBG as $A(t,t_{p})=10$ ($t\le t_{p}$) in original MASs. Meanwhile we set $A(t,t_{p})=1000$ ($t\le t_{p}$) and $\varrho=1$ in MAS (\ref{A1}).
		By Fig. \ref{f7}, we find that, using two MASs, both states converge to a state that is extremely close to the optimal solution at $t_{p}=70$. However, using MAS (\ref{A1}) (CPU time is 6.134 s) is faster than using original MAS (CPU time is 43.812 s) and achieves a solution that is closer to the optimal solution at t=70 (the optimal solution is approximately [2.2090;0.6545;0.1364]). 
		
		\begin{figure}[htph] 
			\centering
			\subfigure[Transient behaviors of $x_{i}(t)$ under original MAS in Example \ref{5e9}. CPU time is 43.812 s. ]{\includegraphics[width=4cm,height=4cm]{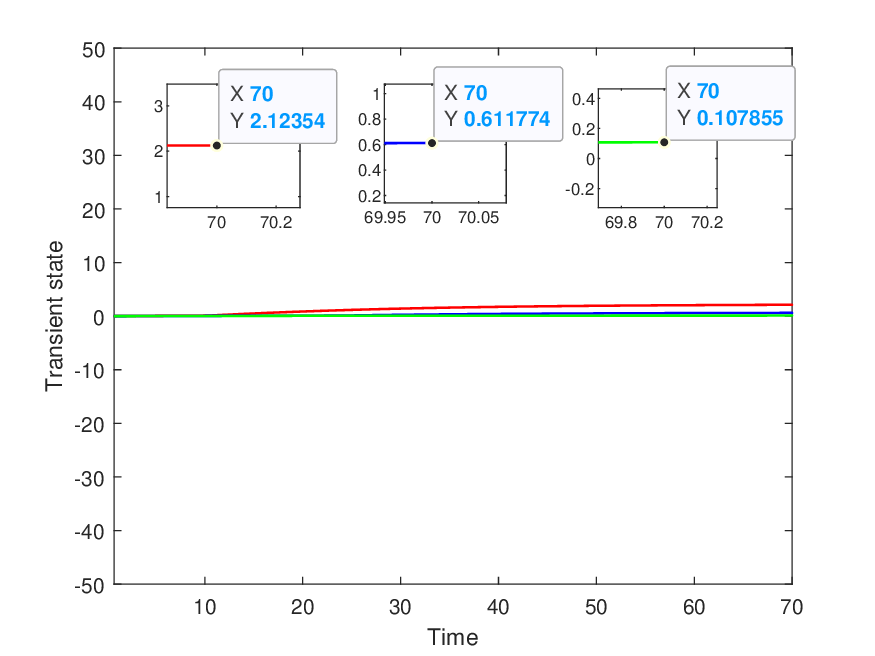}}
			\subfigure[Transient behaviors of $x_{i}(t)$ under MAS (\ref{A1}) in Example \ref{5e9}. CPU time is 6.134s.]{\includegraphics[width=4cm,height=4cm]{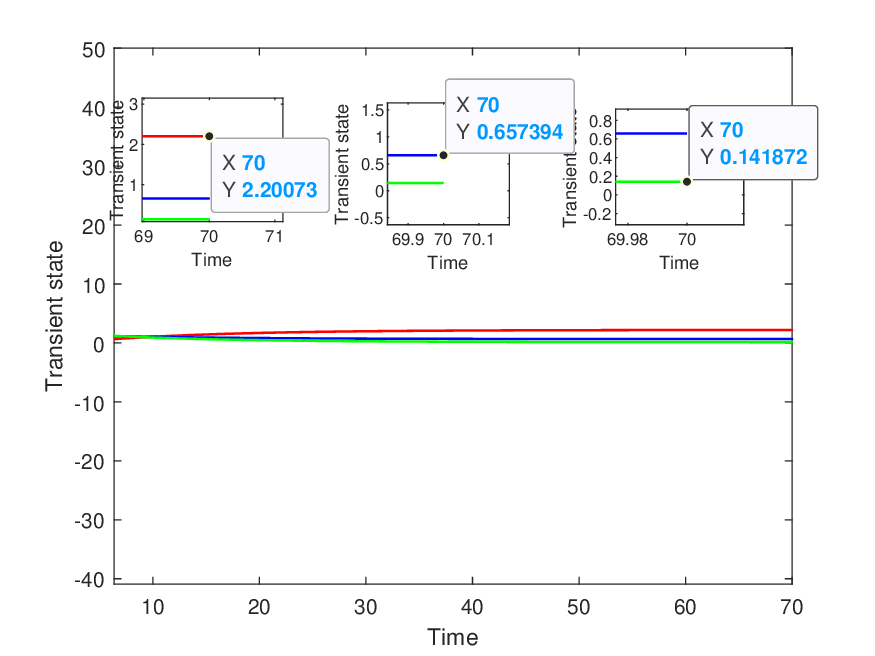}}	
			\subfigure[``old" and ``new" represent the transient behaviors under original MAS and MAS (\ref{A2}) in Example \ref{5e3}, respectively.]{\includegraphics[width=4cm,height=4cm]{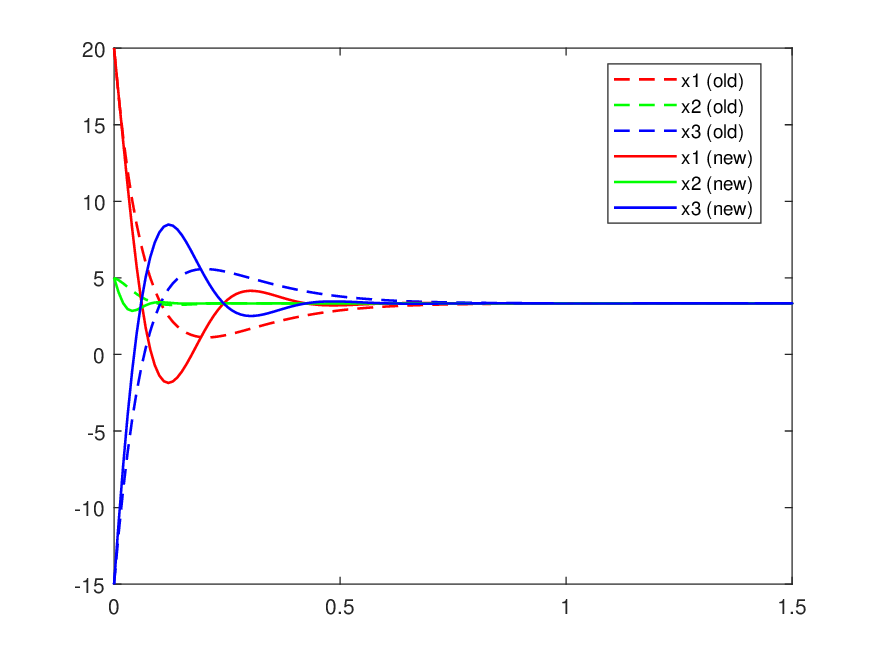}}
			\caption{Transient behaviors of $x_{i}(t)$ in the different MASs}\label{f7}
		\end{figure}
		
	\end{exmp}
	
	\subsection{Examples of consensus-based distributed optimization problems}
	Firstly, we aim to show that under the assumptions stated in (2) of Corollary \ref{c9}, proposed MAS (\ref{A2}) achieves comparable performance to the original MAS in \cite{YL-2023}.
	\begin{exmp}\label{5e3}
		Consider Example 2 in \cite{YL-2023}, and set $f_{1}(x_{1})=\frac{1}{2}x_{1}^{2}-x_{1}$, $f_{2}(x_{2})=\frac{1}{2}x_{2}^{2}-9x_{2}$ and $f_{3}(x_{3})=\frac{1}{2}x_{3}^{2}$.
		Set 
		$x_{1}(0)=20$, $x_{2}(0)=5$, $x_{3}(0) =- 15$,  
		$w_{i}(0)=0$ ($i\in P_{3}$), 
		$L=[1,-1,0;-1,2,-1;0,-1,1]$, and $t_{p}=1$. 
		Now we design the TBG as $A(t,t_{p})=10$ ($t\le t_{p}$) in original MASs. Meanwhile we set $A(t,t_{p})=20$ ($t\le t_{p}$) and $\varrho=0.2$ in MAS (\ref{A2}).
		From Fig. \ref{f7} (c), the states converge to a state that is extremely close to the optimal solution at $t_{p}=1$. Furthermore, Fig. \ref{f7} (c) shows that MAS (\ref{A2}) is also applicable under strong convexity conditions, and the effect does not change much from previous articles.
	\end{exmp}	
	
	Now we aim to show that under the assumptions stated in (3) and (4) of \ref{c9}, the transient state of $x_{i}(t)$ for proposed MAS (\ref{A2}).	
	\begin{exmp}\label{5e10}
		First we consider the examples in Example \ref{5e6}. 
		We set
		$x_{1}(0)=30$, $x_{2}(0)=10$, $x_{3}(0) =-30$,  
		$w_{i}(0)=0$ ($i\in P_{3}$) for Example \ref{5e6},
		Additionally, we set $t_{p}=1.8$ and $L=[1,-1,0;-1,2,-1;0,-1,1]$. 
		Now we design the TBG as $A(t,t_{p})=10$ ($t\le t_{p}$) and $\varrho=0.1$ in MAS (\ref{A1}).
		Fig. \ref{f8} shows that the states converge to a state that is extremely close to the optimal solution at $t_{p}=1.8$.
		Then we consider the following cost functions:
		$\nabla f_{1}(x_{1})=0.8\sin\frac{1}{0.4}-\cos\frac{1}{0.4}(x_{1}\ge0.4)$;
		$\nabla f_{1}(x_{1})=2x\sin\frac{1}{x_{1}}-\cos\frac{1}{x_{1}}(0.2<x_{1}<0.4)$;
		$\nabla f_{1}(x_{1})=0.4\sin\frac{1}{0.2}-\cos\frac{1}{0.2}(x_{1}\le0.2)$.
		$\nabla f_{2}(x_{2})=2\nabla f_{1}(x_{2})$, and $\nabla f_{3}(x_{3})=7\nabla f_{1}(x_{3})$. 
		Set $x_{1}(0)=-50$, $x_{2}(0)=10$, $x_{3}(0) = 50$,  
		$w_{i}(0)=0$ ($i\in P_{3}$). Furthermore, we set $t_{p}=9$ and
		$L=[1,-1,0;-1,2,-1;0,-1,1]$. 
		Now we design the TBG as $A(t,t_{p})=10$ ($t\le t_{p}$)and $\varrho=0.1$ in MAS (\ref{A1}).
		Fig. \ref{f8} shows that the states converge to a state that is extremely close to the optimal solution at $t_{p}=9$.
		
		\begin{figure}[htph] 
			\centering
			\subfigure[Transient behaviors of $x_{i}(t)$ in MAS (\ref{A2}) for Example \ref{5e10}. ]{\includegraphics[width=3cm,height=3cm]{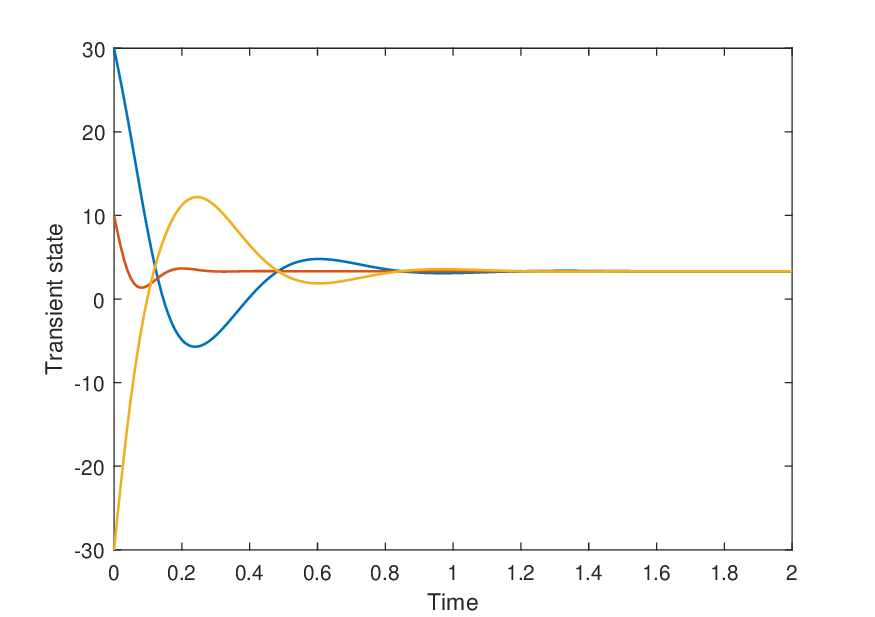}}
			\subfigure[Transient behaviors of $x_{i}(t)$ in MAS (\ref{A2}) for Example \ref{5e10}. ]{\includegraphics[width=3cm,height=3cm]{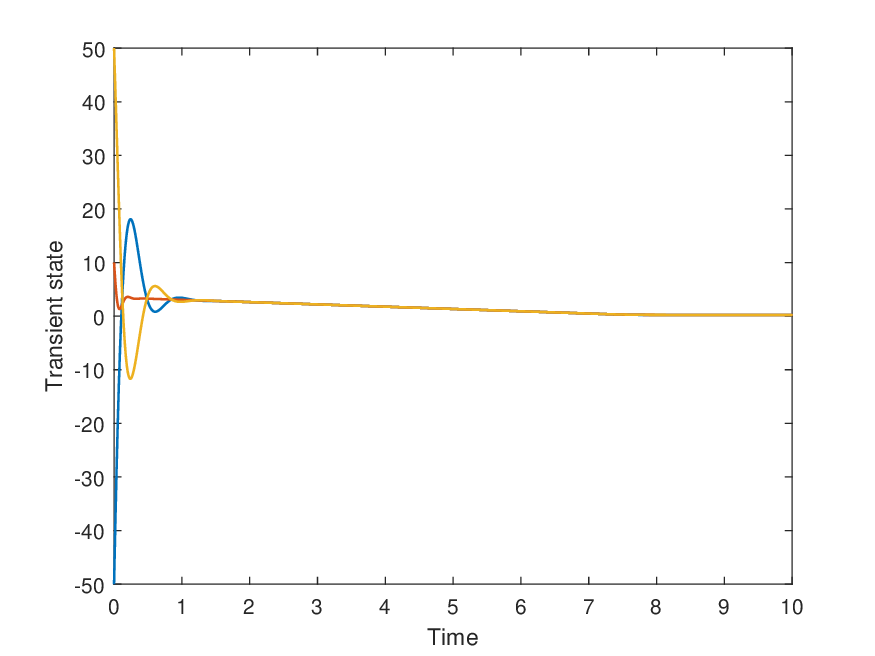}}	
			\subfigure[Transient behaviors of $x_{i}(t)$ in original MAS for Example \ref{5e11}. ]{\includegraphics[width=3cm,height=3cm]{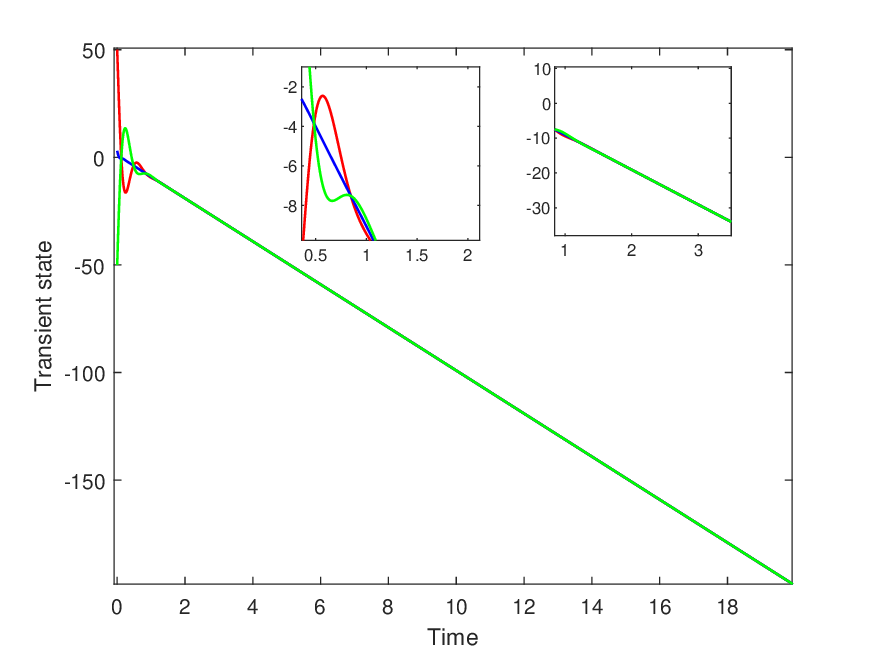}}
			\subfigure[Transient behaviors of $x_{i}(t)$ in MAS (\ref{A2}) for Example \ref{5e11}. ]{\includegraphics[width=3cm,height=3cm]{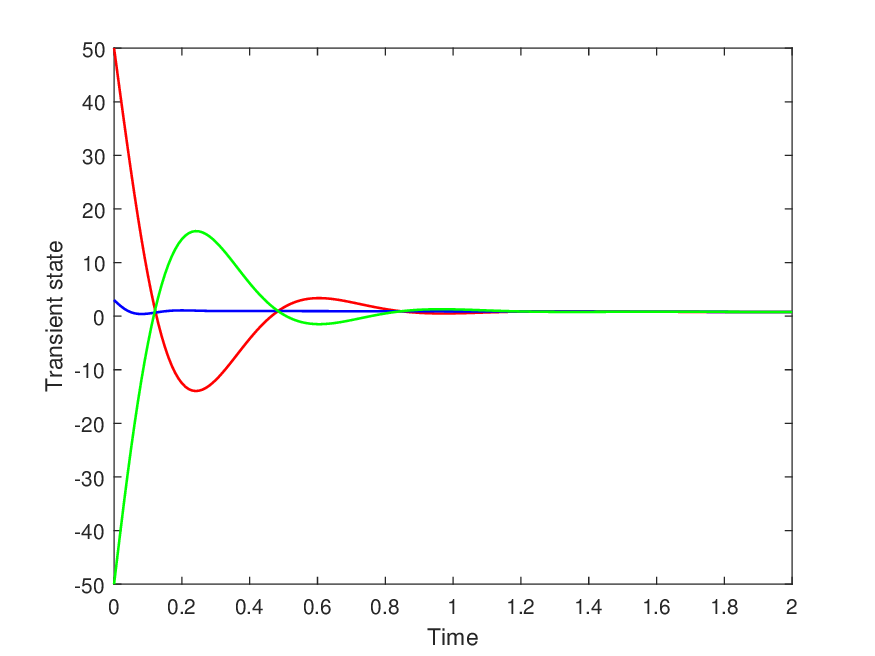}}	
			\caption{Transient behaviors of $x_{i}(t)$ in the different MASs}\label{f8}
		\end{figure}
		
	\end{exmp}
	
	Now we aim to show that under the assumptions stated in (5) of Corollary \ref{c9}, the transient state of $x_{i}(t)$ for MAS (\ref{A2}) and the MAS in \cite{YL-2023}, respectively. 
	
	\begin{exmp}\label{5e11}
		We consider the example:
		$
		f_{1}(x_{1})=x_{1}\sin\frac{1}{x_{1}},
		$
		$f_{2}(x_{2})=f_{1}(x_{2})$ and $f_{3}(x_{3})=f_{1}(x_{3})$.
		We set $x_{1}(0)=50$, $x_{2}(0)=3$, $x_{3}(0) = -50$,  
		$w_{i}(0)=0$ ($i\in P_{3}$). Furthermore, we set $t_{p}=1.8$ and
		$L=[1,-1,0;-1,2,-1;0,-1,1]$.
		Now we design the TBG as $A(t,t_{p})=10$ ($t\le t_{p}$) in original MASs. Meanwhile we set $A(t,t_{p})=10$ ($t\le t_{p}$) and $\varrho=0.1$ in MAS (\ref{A1}).
		Fig. \ref{f8} (c) shows that the variation process of $x$ exhibits an oscillatory pattern and fails to converge. However, Fig. \ref{f8} (d) shows that the states converge to a state that is extremely close to the optimal solution at $t_{p}=1.8$ and the convergence process of $x$ exhibits a relatively smooth trend. Note that, the optimal solution is not unique. Hence, the MAS converges to one of the optimal solutions of the optimization problem.
	\end{exmp}
	
	By replacing $M$-smooth by  smooth stated in (2) of Corollary \ref{c9}, we aim to show the transient state of $x_{i}(t)$ for proposed MAS (\ref{A2}) and the original MAS.
	\begin{exmp}\label{5e13}
		Consider the example in Example \ref{5e9}.	
		We set
		$x_{1}(0)= 50$, $x_{2}(0)= 3$, $x_{3}(0) = -50$,  
		$w_{i}(0)=0$ ($i\in P_{3}$), $t_{p}=2$, $L=[1,-1,0;-1,2,-1;0,-1,1]$.
		Now we design the TBG as $A(t,t_{p})=10$ ($t\le t_{p}$) in original MASs. Meanwhile we set $A(t,t_{p})=10$ ($t\le t_{p}$) and $\varrho=0.1$ in MAS (\ref{A1}).
		By Fig. \ref{f10}, we find that, using two MASs, both states converge to a state that is extremely close to the optimal solution at $t_{p}=2$. However, using MAS (\ref{A2}) (0.525 s) is faster than using original MAS (0.894 s) and achieves a solution that is closer to the optimal solution at $t=2$ (the optimal solution is approximately [0;0;0]). 
		\begin{figure}[htph] 
			\centering
			\subfigure[Transient behaviors of $x_{i}(t)$ in original MAS. CPU time is 0.894 s. ]{\includegraphics[width=4cm,height=4cm]{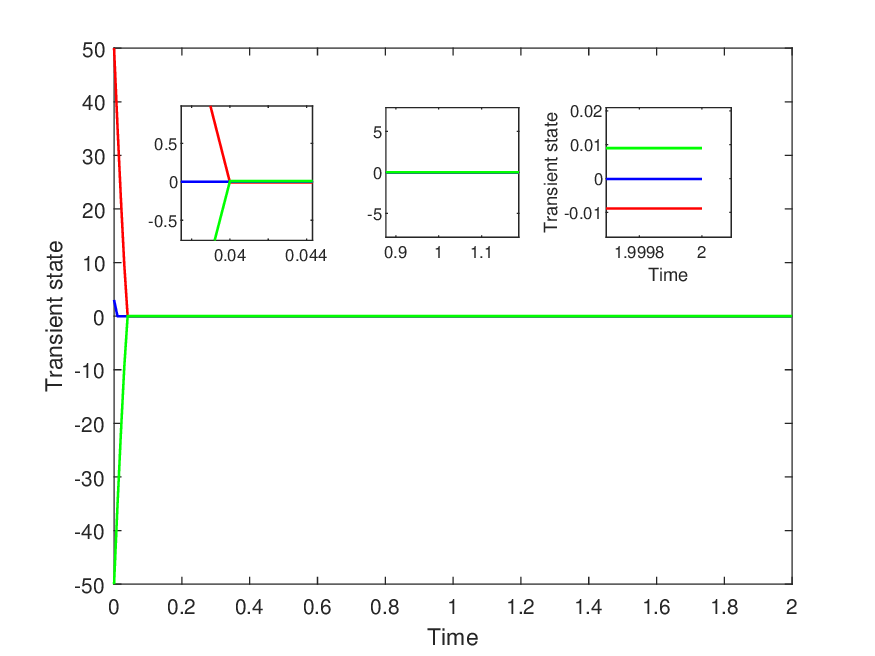}}
			\subfigure[Transient behaviors of $x_{i}(t)$ in MAS (\ref{A2}). CPU time is 0.525 s. ]{\includegraphics[width=4cm,height=4cm]{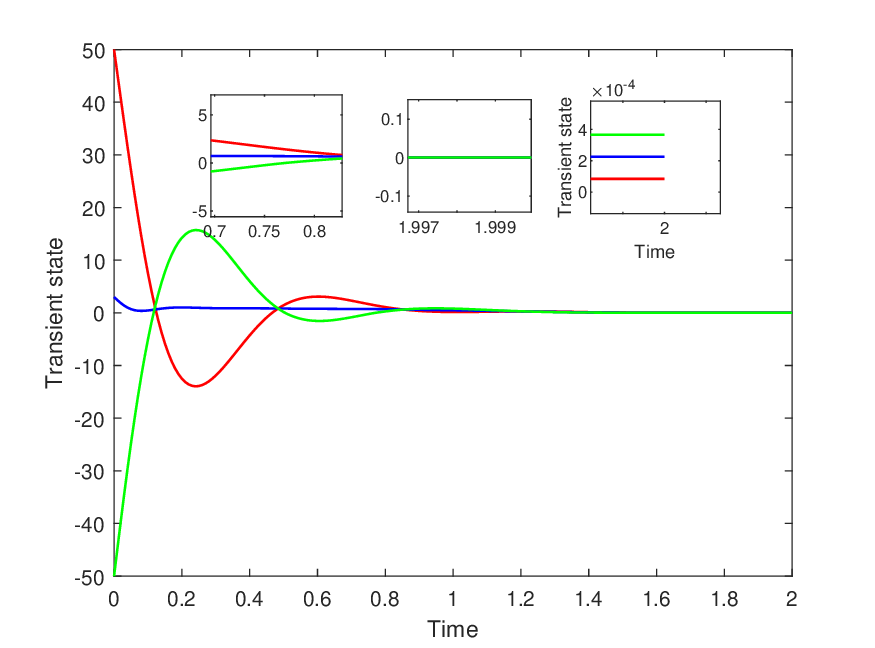}}	
			\caption{Transient behaviors of $x_{i}(t)$ in the different MASs for Example \ref{5e13}.}\label{f10}
		\end{figure}
	\end{exmp}
	
	\section{Conclusion}\label{s6}
	Two MASs under TBG in Definition \ref{aptu} are proposed to solve RAPs and consensus-based distributed non-convex optimization problems within predefined time. Distributed approaches are adopted to enhance security and privacy. The new definitions of predefined-time approximate convergence and predefined-time optimization based on TBG and generalized smoothness are provided. The novel Lyapunov theory for predefined-time approximate convergence is developed. Robustness and boundedness theorems of the Lyapunov function based on TBG in Definition \ref{aptu} are proved. Predefined-time optimization of the proposed MASs is achieved based on the robustness and boundedness theorems. Some numerical examples are presented. Compared with previous MASs, the effectiveness of the calculation results is proved. Future work will focus on addressing global optimal solutions of non-convex functions in distributed optimization.	
	
	\bibliographystyle{siamplain}
	\bibliography{references}

\end{document}